\documentclass[12pt]{article}
\usepackage{amsfonts}
\usepackage{amssymb}
\usepackage{amsthm}
\usepackage[centertags]{amsmath}
\usepackage{newlfont}
\usepackage{graphicx}

\def\r{\mathbb R}

\newtheorem{theorem}{Theorem}[section]
\newtheorem{definition}[theorem]{Definition}
\newtheorem{proposition}[theorem]{Proposition}
\newtheorem{remark}[theorem]{Remark}
\newtheorem{corollary}[theorem]{Corollary}

\title{Invariant surfaces in Sol$_3$ with constant mean curvature and their computer graphics}
\author{Rafael L\'opez\footnote{Partially supported by MEC-FEDER
 grant no. MTM2011-22547 and
Junta de Andaluc\'{\i}a grant no. P09-FQM-5088.}\\
 Departamento de Geometr\'{\i}a y Topolog\'{\i}a\\
Universidad de Granada\\
18071 Granada. Spain\\
email: rcamino@ugr.es}

\date{}

\begin{document}

\maketitle



\begin{abstract} In   Sol$_3$ space there are three uniparametric groups of isometries. In this work we study constant mean curvature  surfaces invariant by one of these groups.  We analyze  the geometric properties of these surfaces by means of their computer graphics. We  construct explicit examples of minimal surfaces and we shall relate them with recent examples of spheres with constant mean curvature.
\end{abstract}

\section{Introduction}

In recent years, the study of surfaces with constant mean curvature in homogeneous three-manifolds is  a topic of great activity, specially after  the extension  by Abresch and Rosenberg  of the Hopf theorem for this kind of spaces (\cite{ar}): see also   \cite{dhm} and references therein. Among the eight models of the geometry of    Thurston (\cite{th}), the space Sol$_3$ is the space with the smallest isometry group.  As a Riemannian manifold,
the space Sol$_3$ can be
represented by $\r^3$ equipped with the metric
$$\langle,\rangle=e^{2z}dx^2+e^{-2z}dy^2+dz^2$$
 where $(x,y,z)$ are the canonical coordinates of $\r^3$. The space Sol$_3$ is a Lie group with the operation
 $$(x,y,z)\ast (x',y',z')=(x+e^{-z}x',y+e^{z}y',z+z')$$
and  the metric $\langle,\rangle$ is  left-invariant with respect to the operator group.

It is natural to consider the study of  curves and surfaces in Sol$_3$ with some added geometric property. For example, the geodesics of Sol$_3$ were studied in \cite{tr},  the totally umbilical invariant surfaces  in \cite{st} and constant angle surfaces in \cite{lm3}. The surfaces in Sol$_3$ with constant mean curvature have received special attention. In the search of explicit examples, some properties on the geometry of the surface  have been assumed, for instance, to be  invariant under   some group of isometries (\cite{lm2}) or  to be a translation surface (\cite{lm1}). Compact surfaces with constant mean curvature and non-empty  boundary were studied in \cite{lo}  and some properties on minimal surfaces in \cite{il}.

Motivated by  the extension of the Hopf theorem, an interesting problem posed in Sol$_3$  was whether exist closed surfaces with constant mean curvature. By the Alexandrov reflection method, an embedded compact surface with constant mean curvature must be of genus $0$. In a first step, Daniel and Mira showed the existence of such surfaces for each value of $H$ with $H>1/\sqrt{3}$ (\cite{dm}). Finally Meeks extended the result for any real number $H$ (\cite{me}).

 The main difficulty   in order to obtain examples of surfaces with constant mean curvature is that Sol$_3$ has a small isometry group as for example, there are no rotations. The isometry group Iso(Sol$_3$) has dimension $3$ and the component of the identity is generated by the following two families of isometries (see \cite{st}):
$$(x,y,z)\longmapsto (\pm e^{-c}x+a,\pm e^c y+b,z+c)$$
$$(x,y,z)\longmapsto (\pm e^{-c}y+a,\pm e^c x+b,-z+c),$$
where $a,b,c\in\r$.

We consider surfaces invariant under a uniparametric group of isometries. In this ambient space, there are three types of such groups, namely, if $\{e_1=(1,0,0),e_2=(0,1,0), e_3=(0,0,1)\}$ is the canonical basis of $\r^3$, each one of the above group is determined by the left translations by $te_i$, $1\leq i\leq 3$:
$${\cal G}_i=\{p\in\mbox{Sol}_3   \longmapsto te_i\ast p; t\in\r\}.$$
Constant mean curvature  surfaces   invariant by ${\cal G}_1$ and ${\cal G}_2$    were classified in \cite{lm2}. In this work we focus in the third group ${\cal G}_3$. Consider
$L_t:\mbox{Sol}_3\rightarrow\mbox{Sol}_3$ the left translation by $te_3=(0,0,t)$:
$$L_t(x,y,z)=(0,0,t)\ast (x,y,z)=(  e^{-t}x, e^t y,t+z).$$
A set $A\subset \mbox{Sol$_3$}$ will called \emph{$3$-invariant} if $L_t(A)\subset A$ for any $t\in\r$. Since we shall  work with immersed surfaces in Sol$_3$, we need to explicit the next definition:

\begin{definition} Let $M$ be a surface and $\psi:M\rightarrow\mbox{Sol$_3$}$ an immersion. We say that $\psi$ a  is $3$-invariant surface if it is the set $\psi(M)$.
\end{definition}

When the immersion is known in the context, we identify $M$ with $\psi(M)$ and we abbreviate by saying that $M$ is a $3$-invariant surface.

In this article we study $3$-invariant surfaces with constant mean curvature $H$. Since a $3$-invariant surface is generated by a curve $\alpha=\alpha(s)$, the condition $H=ct$ expresses as an ordinary differential equation and  for any initial condition we have a solution. However it is difficult to solve the equation $H=ct$  in all its generality, even in the minimal case ($H=0$).  Part of our study is supported in the use of a computer to make numerical pictures of the surfaces that in our case, it has  been possible  by using a symbolic program such as Mathematica. This contrast to with the minimal surfaces invariant by the groups ${\cal G}_1$ and ${\cal G}_2$, where all solutions were done in \cite{lm2}.

The article is organized as follows.    Section \ref{s-2}  takes up the most part of the article where we study $3$-invariant surfaces with zero mean curvature. We shall obtain explicit examples and some  geometric properties. At the end of Section \ref{s-2}, we will relate some of these minimal surfaces   with spheres of Sol$_3$ with constant mean curvature. In Section \ref{s-3}, we consider surfaces with non-zero constant mean curvature.  Due to the difficulty of the mean curvature equation, we   give some graphics of such surfaces, showing by numerical computations the existence of generating curves that are simply closed curves.   Last Section \ref{s-4} is devoted to give   explicit examples of constant Gauss curvature surfaces in Sol$_3$.

\section{Minimal surfaces}\label{s-2}

Let $\psi:M\rightarrow \mbox{Sol}_3$ be a $3$-invariant surface.  The orbit of a point $(x,y,z)$ of Sol$_3$ is the curve $\{L_t(x,y,z); t\in\r\}$, which  acroses once the plane $z=0$ (at the time $t=-z$). Therefore $M$ parametrizes by the set of orbits of a curve $\alpha$ contained in the plane $z=0$. This curve is called the \emph{generating curve} of the surface. Let $\alpha(s)=(x(s),y(s),0)$, $s\in I\subset\r$, be a parametrization of such curve. By the expression of $L_t$, the surface $M$ has a single surface patch $\psi$  given by
\begin{equation}\label{para}
\psi(s,t)=(e^{-t}x(s),e^t y(s),t),\ s\in I, t\in \r.
\end{equation}

\begin{remark}\label{re}
The isometries $\phi (x,y,z)=(\pm e^{-c} y,\pm e^c x,-z+c)$, where $c$ is a real parameter,  carry $3$-invariant surfaces into   new $3$-invariant surfaces. In fact, if  $\alpha$ is the generating curve of a $3$-invariant surface $M$, then $\phi\circ\alpha$ generates the $3$-invariant surface $\phi(M)$. This is because $L_{t}\circ \phi=\phi\circ L_{-t}$.
\end{remark}

We compute the mean curvature of a given $3$-invariant surface $\psi$. First, we recall that  in Sol$_3$ there exists a left-invariant orthonormal frame $\{E_1,E_2,E_3\}$    given by
$$E_1=e^{-z}\frac{\partial}{\partial x},\ \  \ E_2=e^{z}\frac{\partial}{\partial y},\ \  \ E_3=\frac{\partial}{\partial z}.$$
The Riemannian connection $\overset{\sim}{\nabla}$ of Sol$_3$ with respect to this frame is
$$\begin{array}{lll}
\overset{\sim}{\nabla}_{E_1} E_1=-E_3 & \overset{\sim}{\nabla}_{E_1}E_2=0&\overset{\sim}{\nabla}_{E_1}E_3=E_1\\
\overset{\sim}{\nabla}_{E_2} E_1=0 & \overset{\sim}{\nabla}_{E_2}E_2=E_3&\overset{\sim}{\nabla}_{E_2}E_3=-E_2\\
\overset{\sim}{\nabla}_{E_3} E_1=0 & \overset{\sim}{\nabla}_{E_3}E_2=0&\overset{\sim}{\nabla}_{E_3}E_3=0.\\
\end{array}
$$
 Consider $N$ a Gauss map of $M$. The   mean curvature $H$ of $\psi$ with respect to a local parametrization $\psi=\psi(s,t)$ is
 $$H= \frac12\ \frac{eG-2fF+gE}{EG-F^2},$$
where, as usually, $\{E,F,G\}$ and $\{e,f,g\}$ stand for the coefficients of the first and second fundamental form, respectively, that is:
$$E=\langle \psi_s,\psi_s\rangle,\ F=\langle \psi_s,\psi_t\rangle, \ G=\langle \psi_t,\psi_t\rangle.$$
$$e=\langle N,\overset{\sim}{\nabla}_{\psi_s}\psi_s\rangle,\ f=\langle N,\overset{\sim}{\nabla}_{\psi_s}\psi_s\rangle,\
g=\langle N,\overset{\sim}{\nabla}_{\psi_t}\psi_t\rangle.$$
 Without loss of generality, we assume that $\alpha$ is parametrized by the arc length. This means that if $\alpha'=(x',y',0)=x'E_1+y' E_2$, then $x'^2+y'^2=1$. Let $\theta=\theta(s)$ be a differentiable function such that
\begin{equation}\label{e1}
x'(s)=\cos\theta(s),\ y'(s)=\sin\theta(s).
\end{equation}
The first derivatives of the parametrization $\psi$ are:
\begin{eqnarray*}
\psi_s&=&(-e^{-t}x',e^t y',0)=\cos\theta E_1+\sin\theta E_2\\
\psi_t&=&(-xe^{-t},ye^t,1)=-xE_1+yE_2+E_3.
\end{eqnarray*}
The first fundamental form is
 $$E=1,\  F=-x\cos\theta+y\sin\theta, \ G=1+x^2+y^2.$$
Denote $W=EG-F^2$ the determinant of the first fundamental form. Then
 $$W=1+A^2,\ A=x\sin\theta+y\cos\theta.$$
Using the covariant derivatives $\overset{\sim}{\nabla}_{E_i}E_j$, we have
\begin{eqnarray*}\overset{\sim}{\nabla}_{\psi_s}\psi_s&=& -\theta'\sin\theta E_1+\theta'\cos\theta E_2-\cos(2\theta) E_3\\
\overset{\sim}{\nabla}_{\psi_s}\psi_t&=&  (x\cos\theta+y\sin\theta)E_3\\
\overset{\sim}{\nabla}_{\psi_t}\psi_t &=&  -xE_1-y E_2+(y^2-x^2)E_3
\end{eqnarray*}
A unit orthogonal vector field to $M$ is
$$N=\frac{1}{\sqrt{W}}\left(\sin\theta E_1-\cos\theta E_2+A E_3\right).$$
Then the mean curvature of $\psi$ computed with this choice of $N$ is
\begin{equation}\label{mean}
H= \frac{\sin(2\theta)(-x\cos\theta+y\sin\theta)-(1+x^2+y^2)\theta'}{2W^{3/2}}.
\end{equation}

\begin{proposition} The generating curve $\alpha$ of a $3$-invariant minimal surface is obtained by a  solution of the differential equations   \eqref{e1} and
\begin{equation}\label{e2}
\theta'=\sin(2\theta)\frac{-x \cos\theta +y\sin\theta}{1+x^2+y^2}.
\end{equation}
\end{proposition}

The solutions   are obtained provided we give  initial conditions:
\begin{equation}\label{initial}
x(0)=x_0,\ \ y(0)=y_0,\ \ \theta(0)=\theta_0.
\end{equation}
We will denote the solution curve as $\alpha(s;x_0,y_0,\theta_0)$.
A first geometric properties of the generating curve are:

\begin{proposition} Let $\alpha(s;x_0,y_0,\theta_0)$ be a generating curve of a $3$-invariant minimal surface. Then
the domain of $\alpha$ is $\r$. Moreover, a $3$-invariant minimal surface can extend to be complete.
\end{proposition}

\begin{figure}[hbtp]
\begin{center} \includegraphics[width=.4\textwidth]{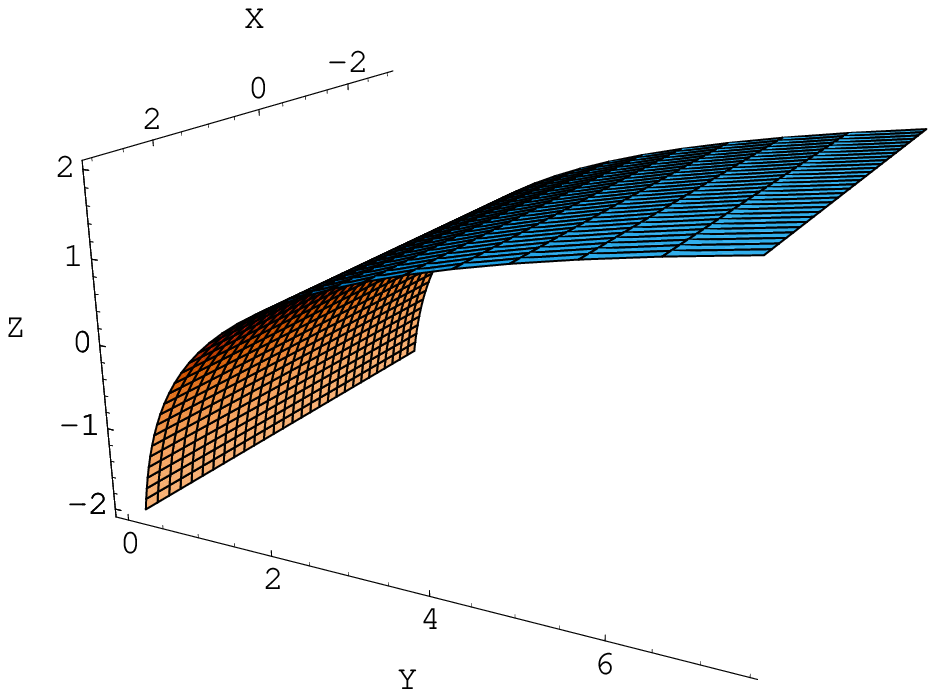}\hspace*{.5cm}\includegraphics[width=.4\textwidth]{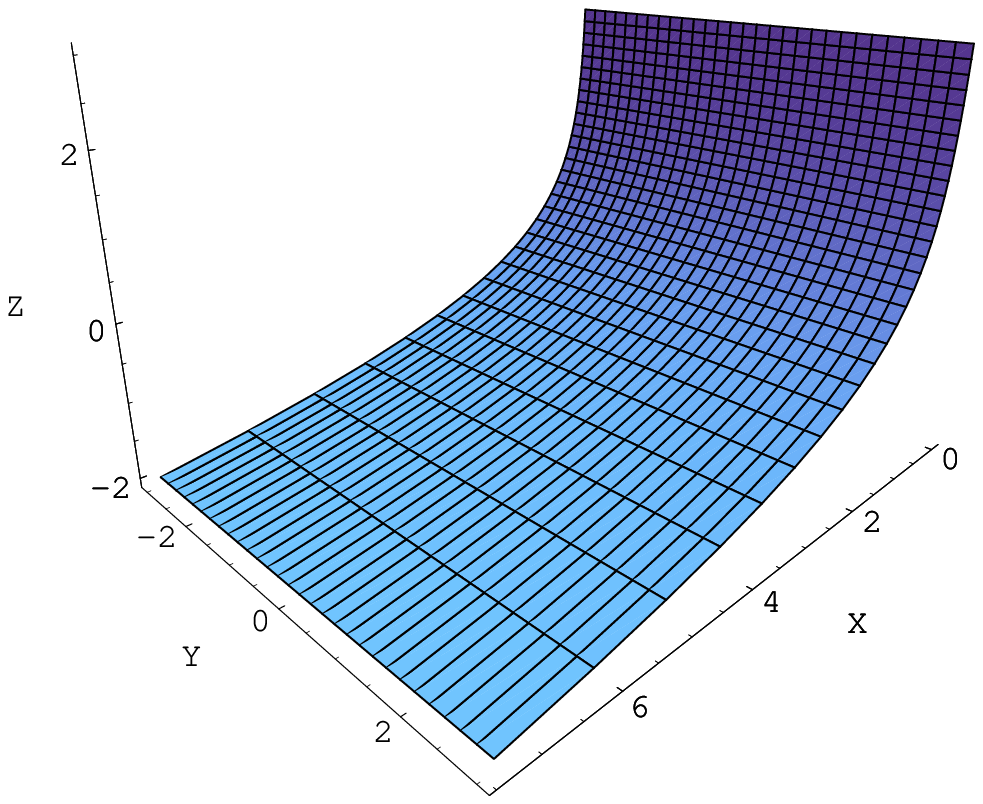}
\end{center}
\caption{(left) A surface of type I parametrized $\psi(u,v)=(u,e^v,v)$; (right) a surface of type II given by $\psi(u,v)=(e^{-u},v,u)$}\label{f1}
\end{figure}

\begin{proof}
The first assertion is a consequence that the derivatives $x'$, $y'$ and $\theta'$ are bounded. For the second statement, we point out that $W=1+A^2\geq 1$. The metric on the surface is the pullback $\psi^{*}(\langle,\rangle)$ of the metric $\langle,\rangle$ of Sol$_3$ by the parametrization $\psi$. As $\psi$ can extend to $\r^2$, we compare $\psi^{*}(\langle,\rangle)$ with  the Euclidean metric $ds_0^2$ of $\r^2$ by
$$ds_0^2=(ds)^2+(dt)^2\leq \psi^*(\langle,\rangle).$$
This proves that the length of any divergent curve is unbounded, in particular, the surface is complete.
\end{proof}
Once obtained the expression of the minimality condition of a $3$-invariant surface, we begin obtaining solutions of \eqref{e2}.   The first examples of such surfaces are the hyperbolic planes $P$ of equation $y=0$ and $Q$  given by
$x=0$. These surfaces are   invariant by the group ${\cal G}_3$ because
$$P={\cal G}_3(\alpha(s)),\ \ \alpha(s)=(s,0,0).$$
$$Q={\cal G}_3(\alpha(s)),\ \ \alpha (s)=(0,s,0).$$
Other  examples of $3$-invariant minimal surfaces appear by choosing $\theta$ a constant function, such it occurs for the planes $P$ and $Q$.
\begin{proposition}\label{pr33}
Given $x_0, y_0\in\r$, we have the next solutions of \eqref{e1}-\eqref{e2}:
\begin{enumerate}
\item  $\alpha(s;x_0,y_0,0)=(s,0,0)+(x_0,y_0,0)$. Here $\theta(s)=0$ (type I).
\item $\alpha(s;x_0,y_0,\pi/2)= s(0,1,0)+(x_0,y_0,0)$. Here  $\theta(s)=\pi/2$ (type II).
\end{enumerate}
Moreover (see Fig. \ref{f1}):
 \begin{enumerate}
\item The surfaces  generated by the curves of type I (resp. type II) are invariant by the group ${\cal G}_1$ (resp. ${\cal G}_2$).
\item Both family of surfaces are ruled surfaces in the sense that they are generated by a uniparametric family of geodesics.
\item There are two   foliations of the  Sol$_3$ space by $3$-invariant minimal surfaces.
\end{enumerate}
\end{proposition}

\begin{proof} We only do the proof for the surfaces of type I. It is immediate that if $\theta_0=0$, then $x(s)=s+x_0$, $y(s)=y_0$ and $\theta(s)=0$ is a solution of \eqref{e1}-\eqref{e2} with initial conditions \eqref{initial}. The parametrization of a such surface is
$\psi(s,t)=(e^{-t}(s+x_0),y_0 e^t,t)$, which can reparametrized as
$$\psi(u,v)=(u,y_0 e^v ,v)=(0,y_0 e^v ,v)+u(1,0,0),\ u,v\in\r.$$
As a consequence, the surface  is invariant by the uniparametric group ${\cal G}_1$ because the isometries of ${\cal G}_1$ are the maps $(x,y,z)\longmapsto (x+t,y,z)$, $t\in\r$. It is also trivial that
 the straight-lines $u\longmapsto \psi(u,v)$ are geodesics of Sol$_3$.
 
Finally,    the foliations by $3$-invariant minimal surfaces are generated by the family of curves  $\{\alpha(s;0,\lambda,0);\lambda\in\r\}$ and $\{\alpha(s;\lambda,0,\pi/2);\lambda\in\r\}$.
\end{proof}

The generating curves of type I (resp. type II) are straight-lines parallel to the $x$-axis (resp. $y$-axis). In particular, the curve $\alpha(s;0,0,0)$ (resp. $\alpha(s;0,0,\pi/2)$) generates the plane $P$ (resp. the plane $Q$). On the other hand, the surfaces of type II can obtain from the ones of type I after  isometries of Sol$_3$, exactly, by $\phi(x,y,z)=(y,x,-z)$.

We recall that all  minimal surfaces in Sol$_3$ invariant by the  uniparametric groups ${\cal G}_1$ and ${\cal G}_2$ were obtained in  \cite{lm2}. Besides the surfaces of type I and II, the classification completes with  the horizontal planes $z=ct$.

We continue obtaining more solutions of  \eqref{e1}-\eqref{e2}. Again, we assume  that $\theta$ is a constant function.

\begin{proposition} Given $x_0\in\r$, the next curves generate $3$-invariant minimal surfaces:
\begin{enumerate}
\item $\alpha(s;x_0,x_0,\pi/4)= s(1,1,0)+(x_0,x_0,0)$, where  $\theta(s)=\pi/4$ (type III).
\item $\alpha(s;x_0,-x_0,-\pi/4)= s(1,-1,0)+(x_0,-x_0,0)$, where  $\theta(s)=-\pi/4$ (type IV).
\end{enumerate}
Moreover both surfaces are ruled, namely,  the curves $s\longmapsto \psi(s,t)$ are geodesic for any $t$.
\end{proposition}
In fact, the curve of type III (resp. type IV) is  the straight-line   $y=x$ (resp. $y=-x$) in the plane $z=0$. The corresponding parametrizations of the surfaces are:
$$\psi(s,t)=(se^{-t},s e^t,t), \ s,t\in\r\ \ \mbox{(type III)}.$$
$$\psi(s,t)=(se^{-t},-s e^t,t), \ s,t\in\r\ \ \mbox{(type IV)}.$$

\begin{proof} We only  consider the surfaces of type III. Again, $x(s)=s+x_0$, $y(s)=s+x_0$ and $\theta(s)=\pi/4$ solve
\eqref{e1}-\eqref{e2} with $\theta_0=\pi/4$. It is  also direct that the coordinate curve $t=ct$ are geodesics  in the ambient space Sol$_3$ hence the surface is ruled.
\end{proof}

\begin{figure}[hbtp]
\begin{center} \includegraphics[width=.8\textwidth]{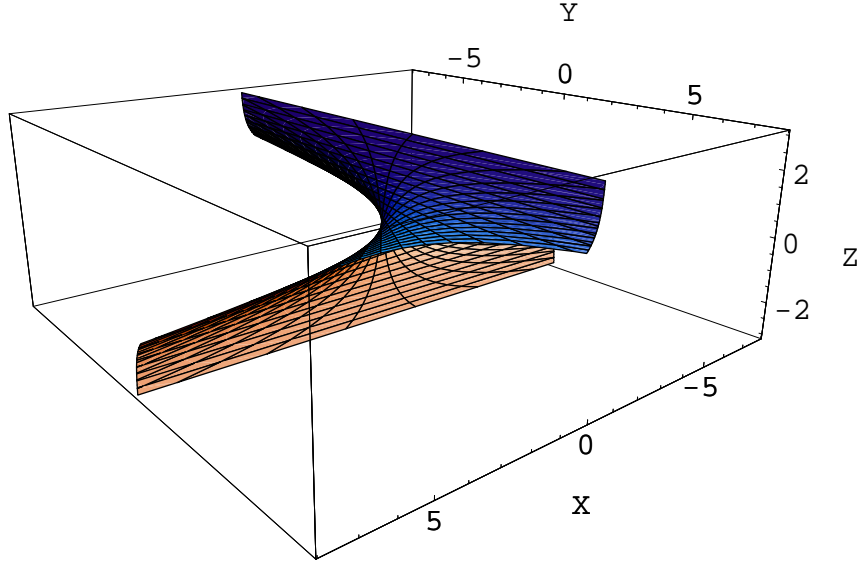}
\end{center}
\caption{The surface of type III parametrized as $\psi(s,t)=(se^{-t},s e^t,t)$.}\label{f2}
\end{figure}

The surfaces of type IV are obtained by the ones of type III after  the isometry $\phi(x,y,z)=(y,x,-z)$. Except for those points of the $z$-axis, the surface of type III (resp. IV) is the Euclidean graph of the function $z=\frac12\log{\left|\frac{y}{x}\right|}$ (resp. $z=\frac12\log{\left|\frac{x}{y}\right|}$): see Fig. \ref{f2}.

\begin{remark} The surfaces of type I, II, III and IV are the only ones that are obtained by choosing $\theta$ a constant function in the solutions of  \eqref{e2}.
\end{remark}

From now, we consider now that $\theta'\not=0$ at some point.

\begin{proposition}\label{pr1} Let $\{x,y,\theta\}$ be a solution of  \eqref{e1}-\eqref{e2}. Then, up $2\pi$ multiple,  the function $\theta$ never attains the values $0,\pi$ or $\pm \pi/2$ unless that $\alpha$ is of type I or type II. In particular, the generating curve is both a graph in the $x$-axis as in  $y$-axis.
\end{proposition}

\begin{proof} We only do the proof in the case that $\theta$ takes the value $0$ at some point $s_0\in\r$.   Define the functions $\{x(s_0)+(s-s_0), y(s_0),0\}$. This  is a solution of  \eqref{e1}-\eqref{e2} such that at $s=s_0$ the  initial conditions are $\{x(s_0),y(s_0),0\}$. Because the same occurs for the functions $\{x,y,\theta\}$, the uniqueness of ODE implies $x(s)=x(s_0)+(s-s_0)$, $y(s)=y(s_0)$ and $\theta(s)=0$ around $s=s_0$. Since this occurs in the maximal interval of definition, we conclude that the curve $\alpha$ is parallel to the $x$-axis, that is, a curve of type I.

Except for the curves of type II and II, the functions $x'(s)$ and $y'(s)$ never vanish in the domain of $\alpha$, which means that $\alpha$ is a graph on the $x$ axis as well as a graph on the $y$-axis.
\end{proof}

 If $\alpha$ is a graph on the axis $x$, $\alpha(x)=(x,y(x),0)$, a direct computation gives   that the minimality condition \eqref{e2} of the corresponding 3-invariant surface generated by $\alpha$ writes as:
\begin{equation}\label{ey}
y''=2y'\frac{yy'-x}{1+x^2+y^2}.
\end{equation}

From this result and with appropriate isometries of Sol$_3$, we can assume, without loss of generality, that the initial velocity vector $\alpha'(0)$ of the generating curve lies in the first quadrant, that is, $\theta(0)=\theta_0\in (0,\pi/2)$. This condition on $\theta_0$ will assumed throughout the rest of the section.

\begin{figure}
 \begin{center} \includegraphics[width=.3\textwidth]{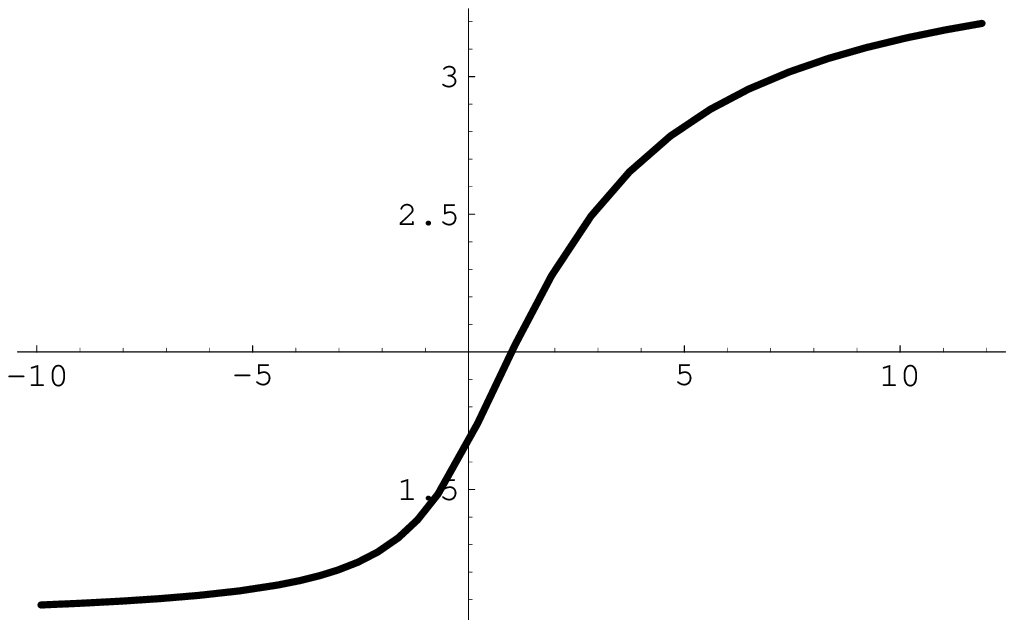}\hspace*{1cm}\includegraphics[width=.3\textwidth]{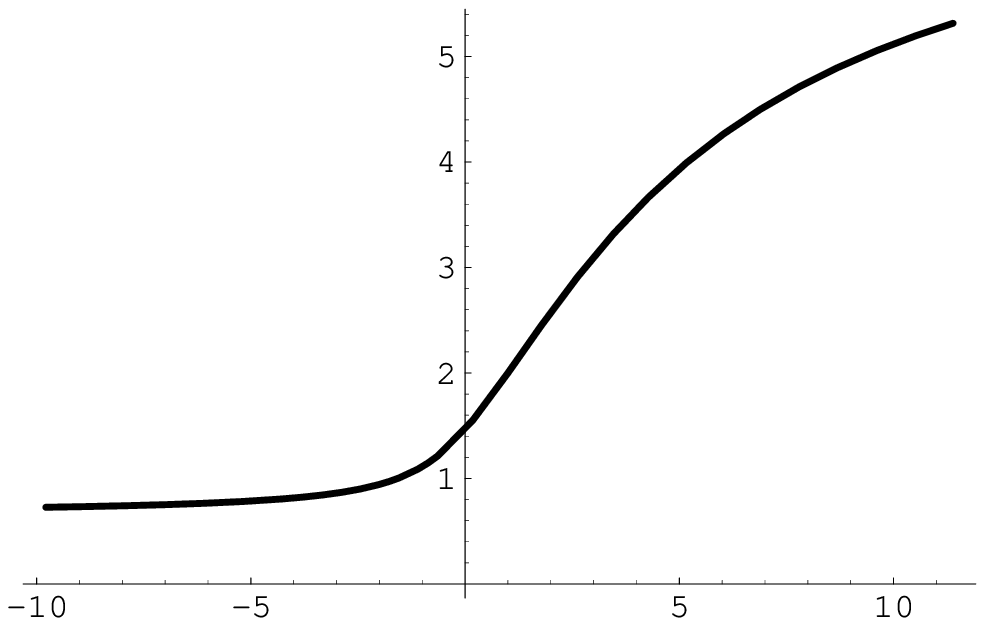}\\
\includegraphics[width=.3\textwidth]{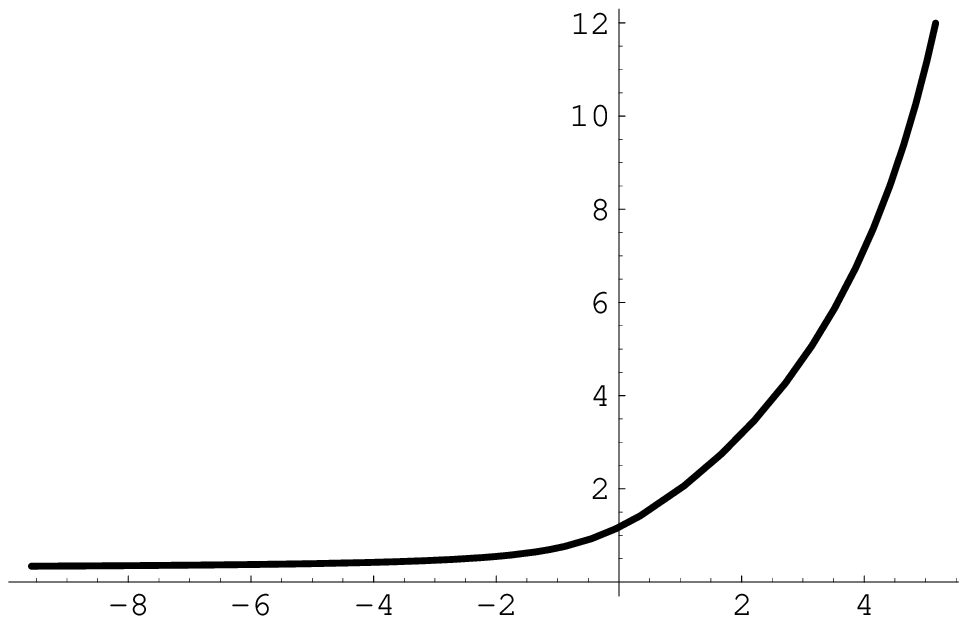}\hspace*{1cm}\includegraphics[width=.3\textwidth]{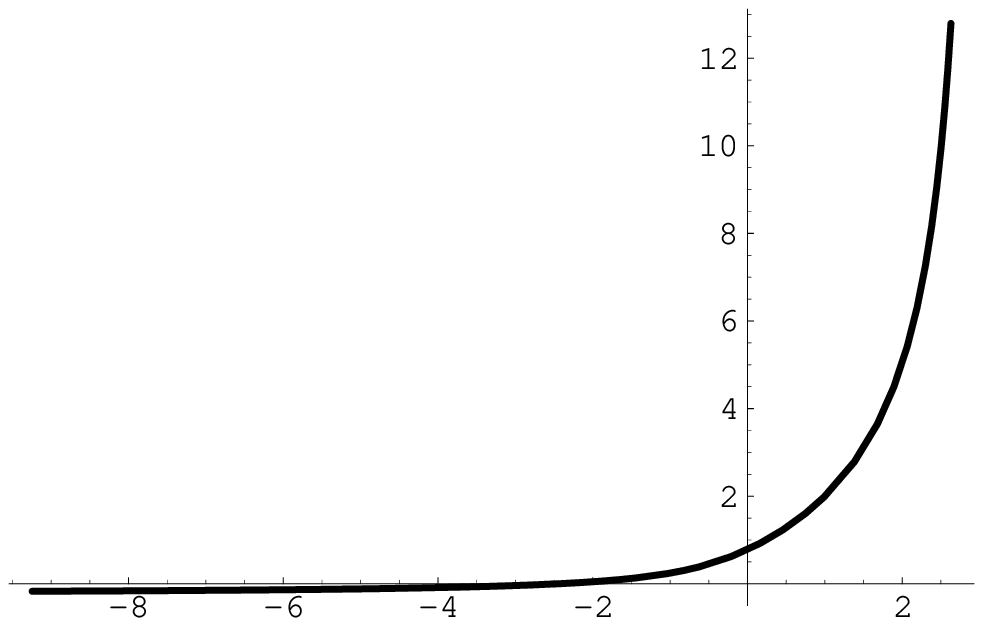}
\end{center}
\caption{We fix the initial values $x_0$ and $y_0$ in \eqref{initial}  and we vary the value of $\theta_0$. Here  $x(0)=1$ and $y(0)=2$. For $\theta_0$, we choose the values  $\pi/10, \pi/6, \pi/4$ and $\pi/3$. At the beginning, the graphic presents one inflection point and it is asymptotic to two parallel lines. After a critical value for  $\theta_0$, the shape of the curve changes to be asymptotic to a quadrant of the $xy$-plane.}\label{f-ex1}
\end{figure}

\begin{figure}
 \begin{center} \includegraphics[width=.3\textwidth]{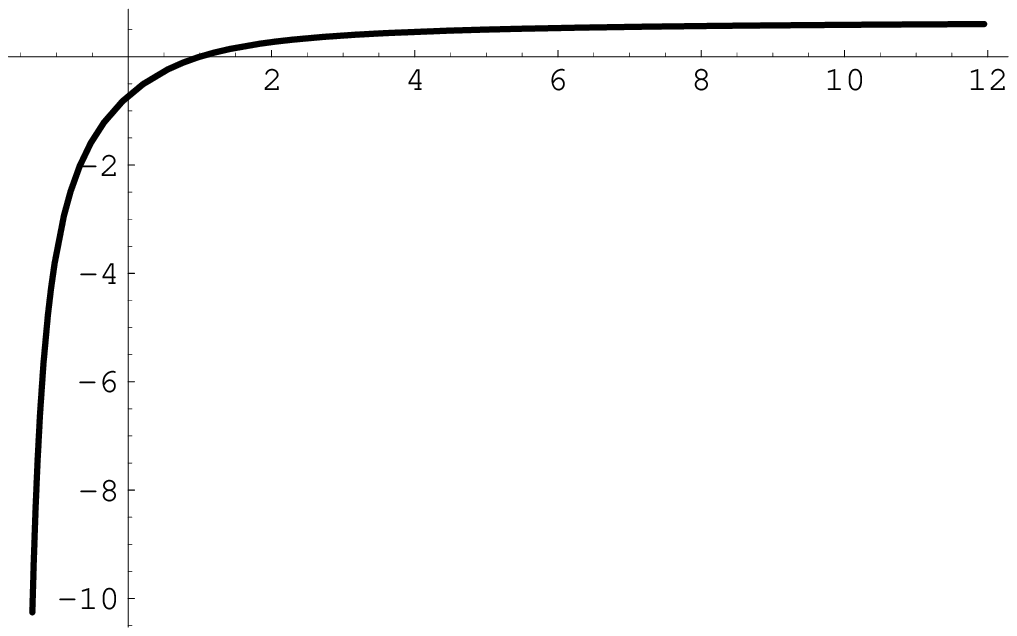}\hspace*{1cm}\includegraphics[width=.3\textwidth]{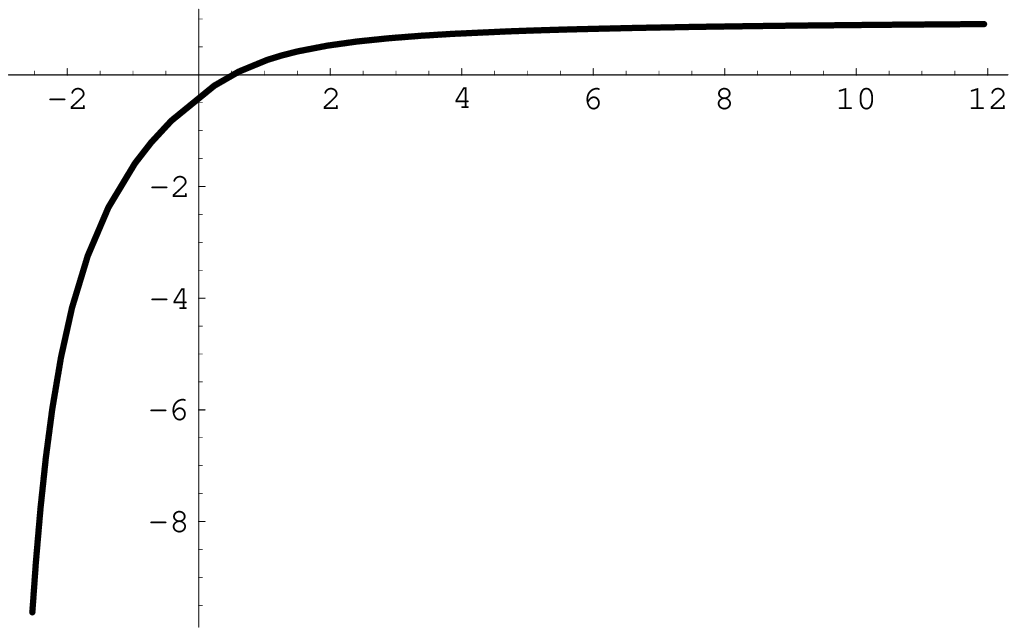}\\
\includegraphics[width=.3\textwidth]{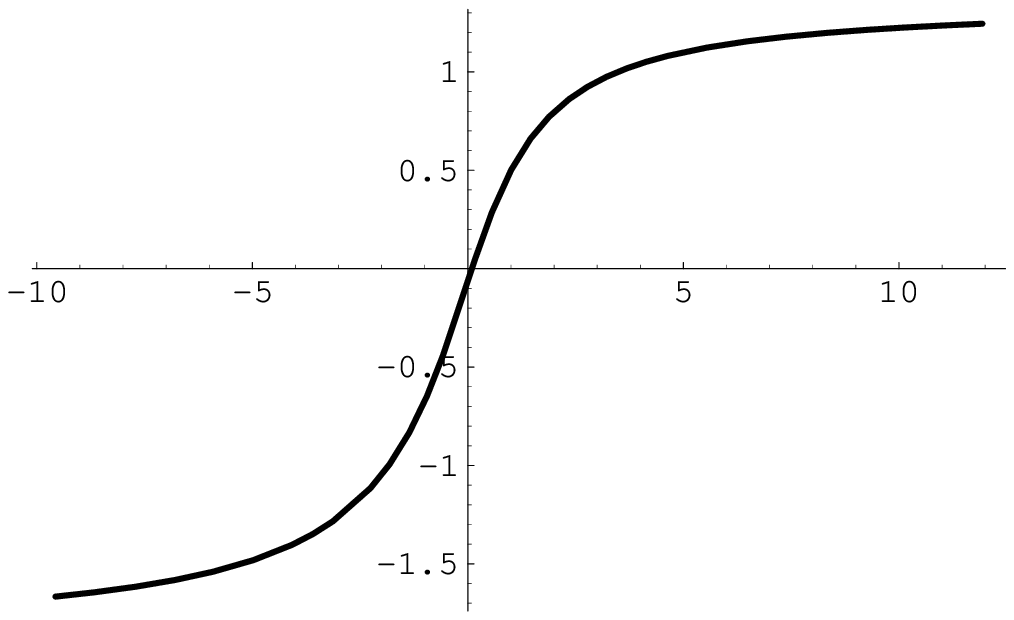}\hspace*{1cm}\includegraphics[width=.3\textwidth]{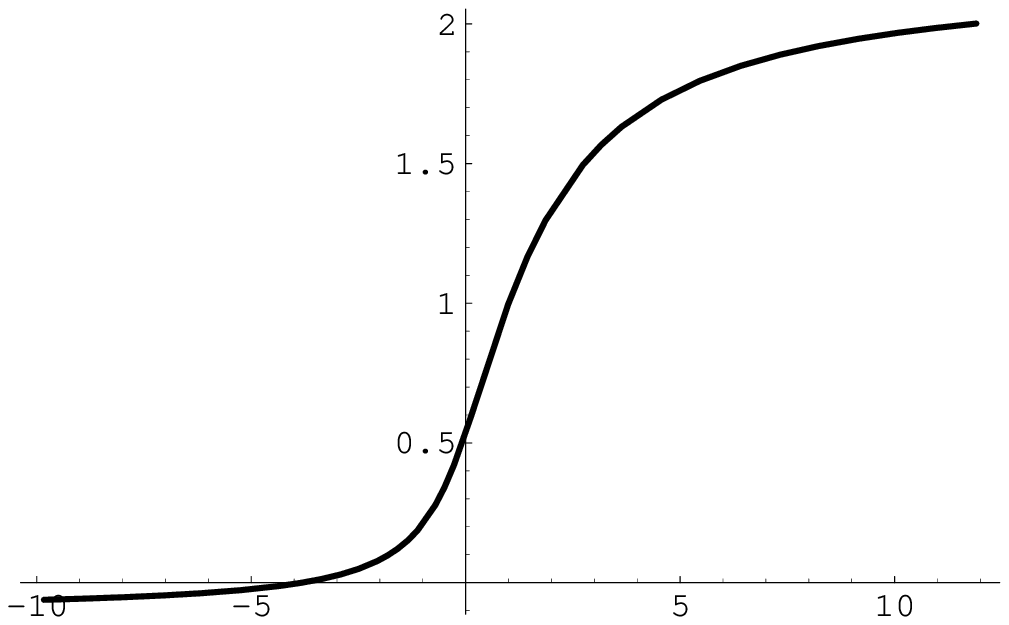}
\end{center}
\caption{In this figure we have fixed  $x_0$ and $\theta_0$ in \eqref{initial}. Here  $x_0=1$ and  $\theta_0=\pi/8$.   The value $y(0)$ takes different values: $0, 1/4, 1/2$ and $1$. In the first two cases, the graphic is asymptotic to a quadrant of the $xy$-plane and for the next two values of $y_0$, the graphic presents one inflection point and it is asymptotic to two horizontal straight-lines.}\label{f-ex2}
\end{figure}

The study of the solutions of \eqref{e2} (or equivalently \eqref{ey}) depends on the initial conditions. We have obtained good numerical pictures of solutions for many initial values. We carried out our graphics using Mathematica and the function {\tt NDSolve} which numerically solves ordinary differential equations. On the basis of our computer pictures, the following result is illustrated in Figs. \ref{f-ex1} and \ref{f-ex2}.

{\bf Experimental result 1}. {\it Depending on the initial solutions and besides the surfaces I--IV, the graphics of solutions of \eqref{e2} are of two types.
\begin{enumerate}
\item (Type A) The curve $\alpha$ is asymptotic to a quadrant of the plane $z=0$ determined by two orthogonal Euclidean straight-lines parallel to the coordinate axis.
\item (Type B) The curve  $\alpha$ is asymptotic to two  straight-lines which are both parallel to the $x$ axis or both parallel to the $y$-axis.
\end{enumerate}
}

In fact, if we fix the initial data $x(0)$ and $y(0)$ and we change the value of $\theta_0$, we expect that there exists a critical value $\theta_0$ where the curve changes from one type to the other one. Assume that the function $y(x)$ describes the generating curve, that is, $y(x)$ is a solution of \eqref{ey}. For the curves of type A, the function $y(x)$ is    convex or concave in all its domain. This would mean that the second derivative $y''(x)$ never vanish. On the other hand, the curves of type B  would have only one inflection point, that is, the equation $y''(x)=0$ would have only one solution $x_0$ and in each one of the intervals $(-\infty,x_0)$ and $(x_0,\infty)$ the curve would change from convex (or concave) shape to concave (or convex).

In what follows, we assume that the   starting point of the generating curve $\alpha$ is the origin $(0,0)$  of $\r^2$ with initial velocity $\theta(0)=\theta_0\in\r$, that, is,
\begin{equation}\label{co}
x(0)=0,\ y(0)=0, \theta'(0)=\theta_0.
\end{equation}
In particular, the $z$-axis is contained in the surface because it is the curve ${\cal G}_3(0,0,0)$.

\begin{theorem} Let $\alpha$ be a generating curve of a $3$-invariant minimal surface satisfying  the initial values \eqref{co}.
\begin{enumerate}
\item The curve $\alpha$ is symmetric with respect to the origin $(0,0)$.
\item The symmetry of $\alpha$ with respect to the line $y=x$ is a generating curve of other $3$-invariant minimal surface.
    \end{enumerate}
\end{theorem}

\begin{figure}[hbtp]
\begin{center} \includegraphics[width=.4\textwidth]{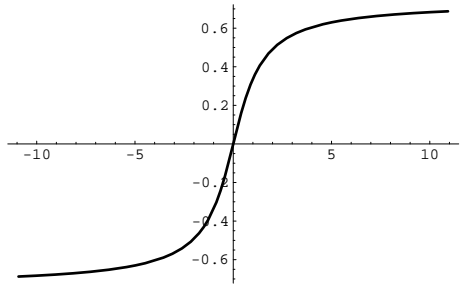}
\end{center}
\caption{A generating curve for initial conditions $x(0)=y(0)=0$ and $\theta(0)=\pi/8$. This curve is asymptotic to the lines $y\simeq \pm 0.736872$.}\label{f-B}
\end{figure}

\begin{proof} Denote $\alpha(s)=(x(s),y(s),0)$ a solution of \eqref{e1}-\eqref{e2} with initial conditions \eqref{co}.
\begin{enumerate}
\item The functions $\bar{x}(s)=-x(-s)$, $\bar{y}(s)=-y(-s)$ and $\bar{\theta}(s)=\theta(-s)$ are a solution of \eqref{e1}-\eqref{e2} with the same initial conditions \eqref{co}.  Let us apply uniqueness of ODE.
\item If $\theta_0=\pi/4$, we know that $\alpha$ is the curve of type III, which agrees with its symmetry with respect to the line $y=x$. We assume $\theta_0\in (0,\pi/4)$. The symmetry of $\alpha$ is the curve $\phi\circ\alpha$, where $\phi(x,y,z)=(y,x,-z)$. Now we use Remark \ref{re}.
\end{enumerate}
\end{proof}
Because the curve is symmetric with respect to the origin, then the curve is of type B depicted above, that is,
 the generating  curve is asymptotic to two parallel straight-lines: see Fig. \ref{f-B}.

By the second statement of this result, we can assume that $\theta_0\in (0,\pi/4)$.

\begin{theorem}\label{t-minimal} Assume $\theta_0\in (0,\pi/4)$ and let $y=y(x)$ be the function  which gives a generating curve $\alpha$ with initial conditions \eqref{co}.
\begin{enumerate}
\item The function $y(x)$ is increasing in all its domain.
\item There exists only one inflection point, namely, $x=0$. If $x>0$ (resp. $x<0$), then $y$ is concave (resp. convex).
\item The curve $\alpha$ does not across the line $y=x$.
\end{enumerate}
\end{theorem}

\begin{proof} By the symmetry of $\alpha$ with respect to the origin, it suffices to do the study of the function $y(x)$ for $x\geq 0$. As $\theta_0\in (0,\pi/4)$, then $0<y'(0)<1$. In particular, $y'>0$ in some interval $(0,\epsilon)$, $\epsilon>0$. If at some point, the derivative $y'$ vanishes, Prop. \ref{pr1} asserts that $y$ is a constant function, which is impossible since $y'(0)\not=0$. As a conclusion, $y'\not=0$ and so, $y'>0$ in $(0,\infty)$.

On the other hand, we know that $y''(0)=0$. Define $h(x)=y(x)y'(x)-x$, which, by \eqref{ey}, gives us the sign of $y''$. Then $h(0)=0$ and $h'(0)=y'(0)^2-1<0$. This means that $h$ is decreasing in some interval $(0,\delta)$, $\delta>0$. From \eqref{ey}, $y''<0$ in  $(0,\delta)$.

{\it Claim.} The function  $h$ is negative in $(0,\infty)$.  On the contrary, let $x_1>0$ be the first point such that $h(x_1)=0$. Then there exists $x_2\in (0,x_1)$ such that $h'(x_2)=0$, that is,
$$y''(x_2)=\frac{1-y'(x_2)^2}{y(x_2)}.$$
As $y'$ is decreasing on $(0,x_1)$, $y'(x_2)<1$ and so, $y''(x_2)>0$. This is a contradiction because $y''$ is negative in $(0,x_1)$. This proves the claim.
 
 Then $y''<0$ in $(0,\infty)$ and $y$ is concave.  As a consequence, the function $y'$ is decreasing in $x$: $y'(x)<y'(0)<1$ on $(0,\infty)$. Near to $x=0$ and for positive values of $x$, the curve $\alpha$  lies  in the domain determined by the lines $y=x$ and $y=0$. If the curve $\alpha$ attains the line $y=x$, by  the intermediate value theorem, $y'=1$ at some point of the interval $(0,\infty)$, which is impossible.
\end{proof}

Recall that the preliminary numerical graphics of the solutions  with initial  conditions \eqref{co} suggest the following result:

{\bf Experimental result 2}. {\it Assume that $\alpha$ is a solution of \eqref{e1}-\eqref{e2} with initial conditions \eqref{co} and  $\theta(0)\in (0,\pi/4)$. Then  the curve $\alpha$ lies between two parallel straight-lines to the $x$-axis, being $\alpha$ asymptotic to both ones.
}

If the initial angle $\theta_0$ moves in the range $(\pi/4,\pi/2)$, then $\alpha$ is asymptotic to two straight-lines parallel to the $y$-line.

Once assumed this result, let us consider $\alpha$ the generating curve of a $3$-invariant minimal surface and  let $L_1$ and $L_2$ be the two straight-lines in the plane $z=0$ which $\alpha$ is asymptotic to. If we apply the group ${\cal G}_3$ to $L_1\cup L_2$, we obtain two $3$-invariant surfaces $P_i={\cal G}_3(L_i)$, $i=1,2$, which lie separated by a constant distance. Then $P_1\cup P_2$ determine a {\it slab} in Sol$_3$. Because $\alpha$ lies between the strip of the plane $z=0$ bounded by $L_1$ and $L_2$, the surface that generates lies in the slab defined by $P_1$ and $P_2$. Thus we conclude

\begin{corollary} In Sol$_3$ there exist many minimal surfaces contained in a slab.
\end{corollary}

This result is the version in Sol$_3$ of  the classical result of Jorge and Xavier about the existence of complete minimal surfaces in a slab  in Euclidean space (\cite{jx})

We finish this section relating the minimal surfaces of type III and IV with the spheres of constant mean curvature and whose existence was proved in \cite{dm,me}: see also \cite{mmpr} for any metric in Sol$_3$. For each $H>0$ there exists (up left translations) a unique sphere with constant mean curvature $H$ which  is embedded. We abbreviate by saying a $H$-sphere. Each one of these $H$-spheres admits certain symmetries of order $2$. Exactly, at each point $p\in \mbox{Sol}_3$  there exists an orthonormal basis $\{v_1(p),v_2(p),v_3(p)\}$ such that the Ricci tensor diagonalizes. It is known that the Gauss map $N$ of a $H$-sphere is a diffeomorphism on the unit sphere. Therefore in a given $H$-sphere there exists points where the Gauss map is a vector that diagonalizes the Ricci tensor, that is, $N(p)=\pm v_i(p)$, for some $1\leq i\leq 3$. In fact, in the $H$-sphere there exist exactly $6$ points. It is proved that  the sphere is symmetric with respect to the rotation by  angle $\pi$ around the normal at $p$ (see also \cite{mmpr}). At the origin $(0,0,0)$, the basis that diagonalizes the Ricci tensor is $\{(1/\sqrt{2})(E_1+E_2), (1/\sqrt{2})(E_1-E_2), E_3\}$.

Take $H_n$ a sequence of positive numbers converging to $0$ and the corresponding $H_n$-spheres $M_n$. For each $M_n$, there exists a point $p_n$ where the unit normal vector is the translated of $(1/\sqrt{2})(E_1-E_2)$. Let us do left translations of each $M_n$ such that $p_n$ agrees with the origin
$p=(0,0,0)$ of $\r^3$ and thus $N_n(p)=(1/\sqrt{2})(E_1-E_2)$. Then all   spheres $M_n$ are invariant with respect to the rotation by angle $\pi$ around the normal at $p$ that fixes $p$. In our case, the tangent plane $T_pM_n$ at $p$ is spanned by
$\{(1/\sqrt{2})(E_1+E_2),  E_3\}$ and the isometry that fixes $N_n(p)$ and which is a rotation by angle $\pi$ in $T_pM_n$ is  given by  $\phi(x,y,z)=(-y,-x,-z)$. As $H_n\rightarrow 0$, the sequence $M_n$ converges to a complete stable minimal surface $M$.  Therefore, the limit surface, which is a minimal surface,   inherits the same symmetry, that is, it is invariant by the action of $\phi$. Of course, the normal vector of this surface at $p$ is $(1/\sqrt{2})(E_1-E_2)$.

On the other hand, the surface of type III whose parametrization is $\psi(s,t)=(se^{-t},se^t,t)$, $s,t\in\r$ has the same tangent plane and normal vector as $M_n$ at the point $p$. Moreover, this surface is invariant by the isometry $\phi$ because
$$\phi(\psi(s,t))=(-se^t,-se^{-t},-t)=\psi(-s,-t).$$
A similar process can do with  the surfaces of type IV. As conclusion of the above reasoning,

\begin{corollary} The surfaces of type III and IV  are  candidate minimal surfaces to be the limit of a sequence of closed $H_n$-spheres, with $H_n\rightarrow 0$.
\end{corollary}

Therefore it is an open question if there exists a sequence of $H_n$-spheres $M_n$, $H_n\rightarrow 0$, such that $M_n$ converges to a surface of type III or IV.

\section{Computer graphics of surfaces with constant mean curvature}\label{s-3}

\begin{figure}[hbtp]
\begin{center} \includegraphics[width=.4\textwidth]{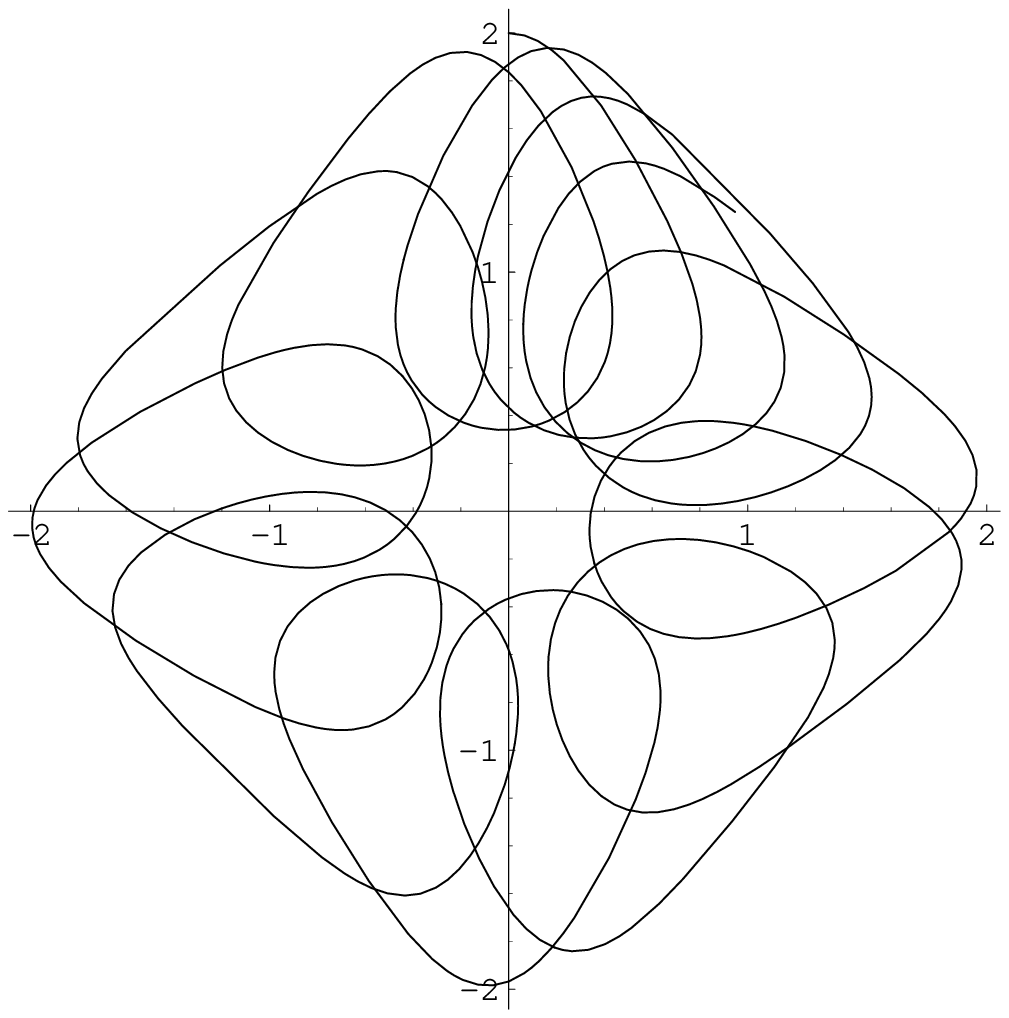}\hspace*{.5cm}\includegraphics[width=.4\textwidth]{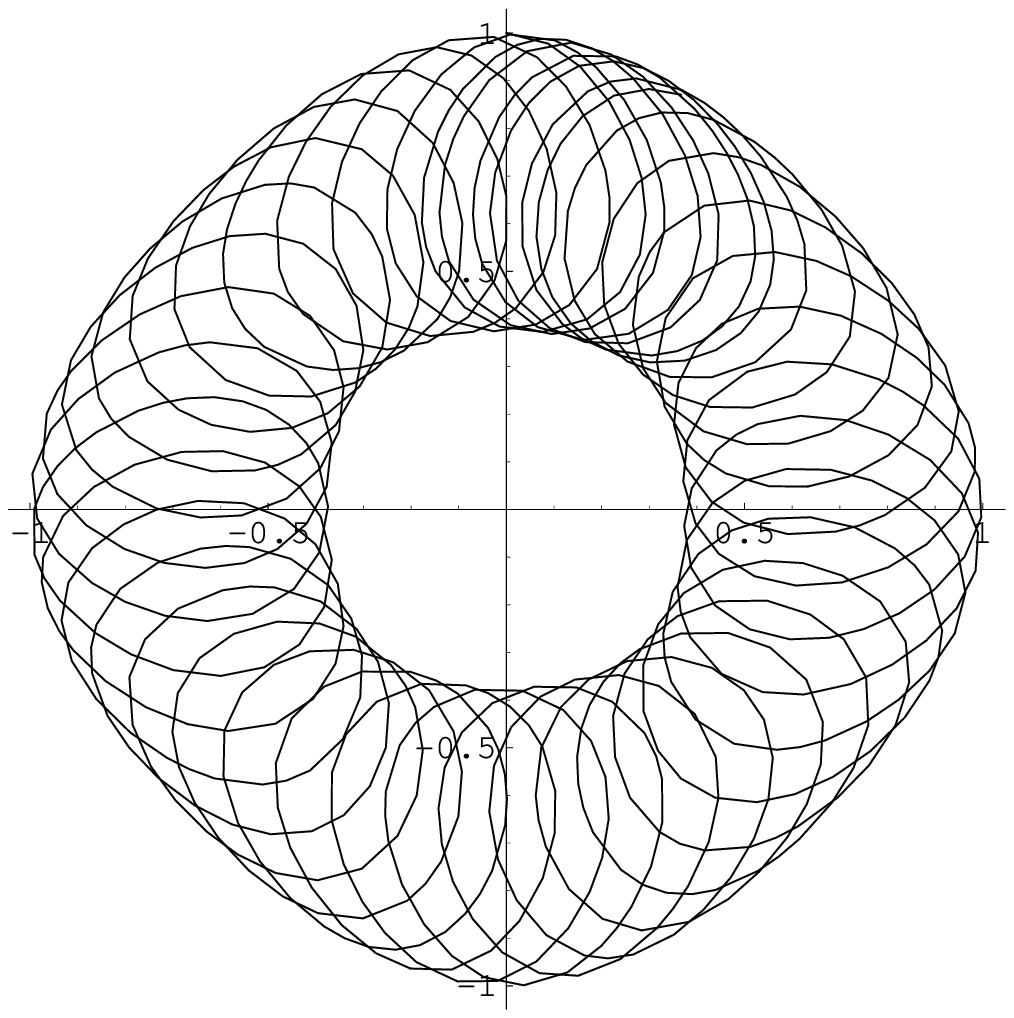}
\end{center}
\caption{Two generating curves of $3$-invariant surfaces with constant mean curvature: (left) $H=1$ and $\alpha(0)=(0,2,0)$; (right) $H=2$ and $\alpha=(0,1,0)$.}\label{f-h1}
\end{figure}

We now consider $3$-invariant surfaces with constant mean curvature $H$ and $H\not=0$. By using the value of $H$ in \eqref{mean} and that $W=1+A^2$, the equation that satisfies the angle function $\theta$ of the generating curve $\alpha$ is
$$\theta'=\frac{1}{1+x^2+y^2}\Big(\sin(2\theta)(-x\cos\theta+y\sin\theta)- 2H(1+A^2)^{3/2}\Big).$$
In this case, the term $2H(1+A^2)^{3/2}$ adds extra difficulties in the study of the geometric properties of the solution curves. By means of the computer, one can produce approximate solutions. In Fig.  \ref{f-h1}, we present two solutions for different initial data. In both cases, the initial angle $\theta_0$ is the same, namely, $\theta_0=0$.

Although the variety of cases, our numerical computations give evidences, as can be seen in the the accompanying computer graphics. What we observe is the following. We fix the value of the (constant) mean curvature, namely, $H=1$ and the initial velocity of the curve at $s=0$, $\theta_0=0$, that is, $\alpha'(0)$ is the horizontal vector $(1,0)$ viewing the curve $\alpha$ in the $xy$-plane. With respect to the initial conditions \eqref{initial}, we take $x(0)=0$ and $y(0)=y_0$, where $y_0$ is a one real parameter. This means that $\alpha(0)$ belongs to the $y$-axis.

The twelve curves pictured below were produced with different values of $y_0$. The rows correspond to the values
$1/8, 1/4, 1/2$ and $3/4$ of the parameter $y_0$ and the columns correspond to the values $2,4$ and $6$ of $s_0$, where $[0,s_0]$ is the domain of the generating curve $\alpha=\alpha(s)$. We observe that as $s_0$ increases, the trace of $\alpha$ begins doing a turn in such way that $\alpha$ acroses the positive $y$-axis again, and many times as $s_0$ goes to $\infty$. We have pointed in the graphic with a black small square the initial point of $\alpha$, that is, $\alpha(0)$, and with a black small disc the following  point where $\alpha$ has tangent horizontal vector $(1,0)$. 
We denote  this new point  as $\alpha(s_1;y_0)=(x(s_1;y_0),y(s_1;y_0))$ where $s_1$ indicates the time where it reaches this point and $y_0$ the dependence on $y_0$. The corresponding velocity vectors at $s=0$ and $s=s_1$ appear as a non-dashed and dashed arrow, respectively. For values $s_0$ close to $0$, $x(s_1;y_0)<0$ and $y(s_1;y_0)>y_0$. However, as $y_0$ increases, the function $x(s_1;y_0)$ is increasing on $y_0$ and $y(s_1;y_0)$ is decreasing on $y_0$. The pictures indicate that we attain a value $y_0^*$ such that
$$\left\{\begin{array}{ll}
& \alpha(s_1;y_0^*)=\alpha(0;y_0^*)\\
& \alpha'(s_1;y_0^*)=\alpha'(0;y_0^*).
\end{array}
\right.$$

 This would say that there exists a \emph{closed simply generating curve}.

 \begin{center} \includegraphics[width=.15\textwidth]{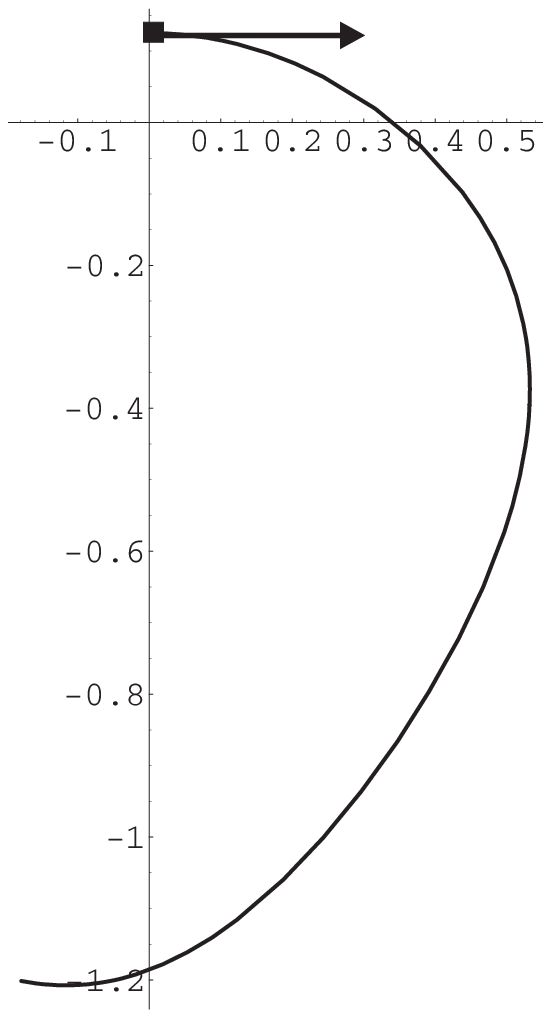}\hspace*{.5cm}\includegraphics[width=.3\textwidth]{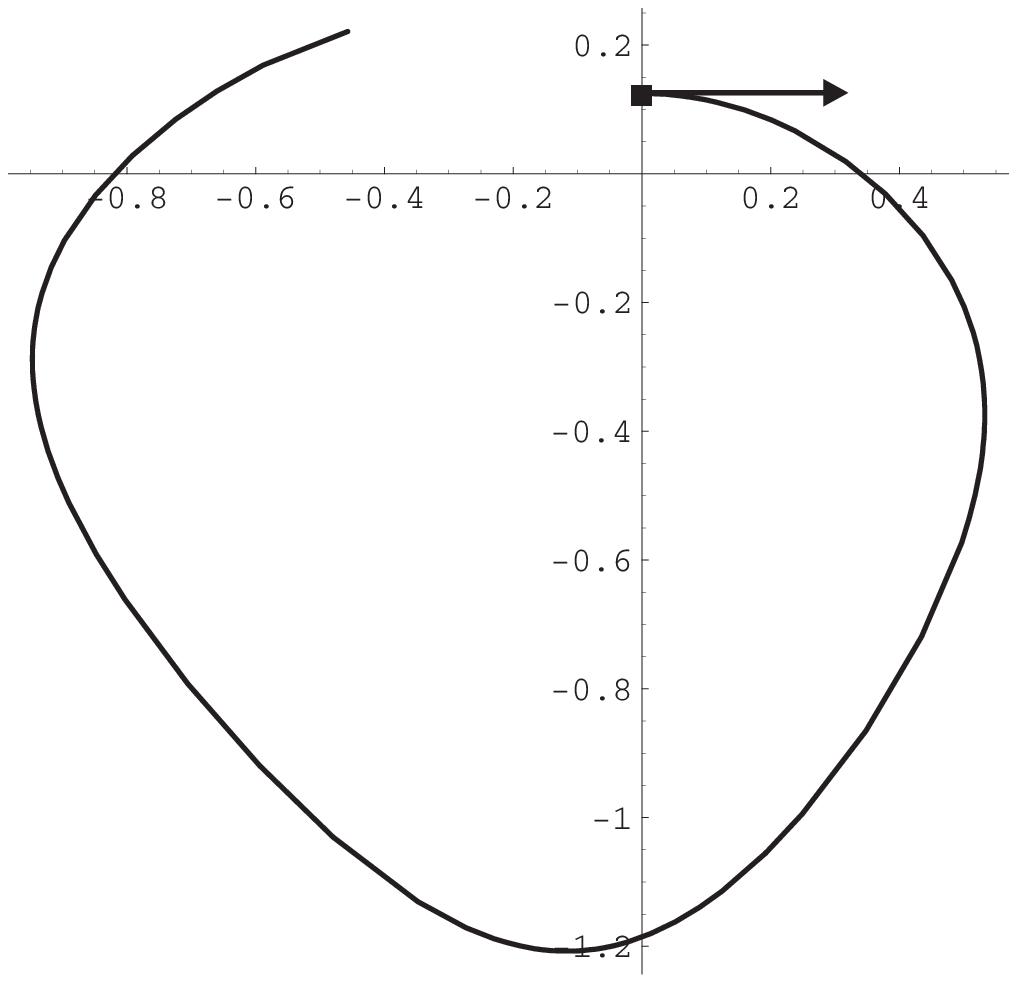}\hspace*{.5cm}
\includegraphics[width=.3\textwidth]{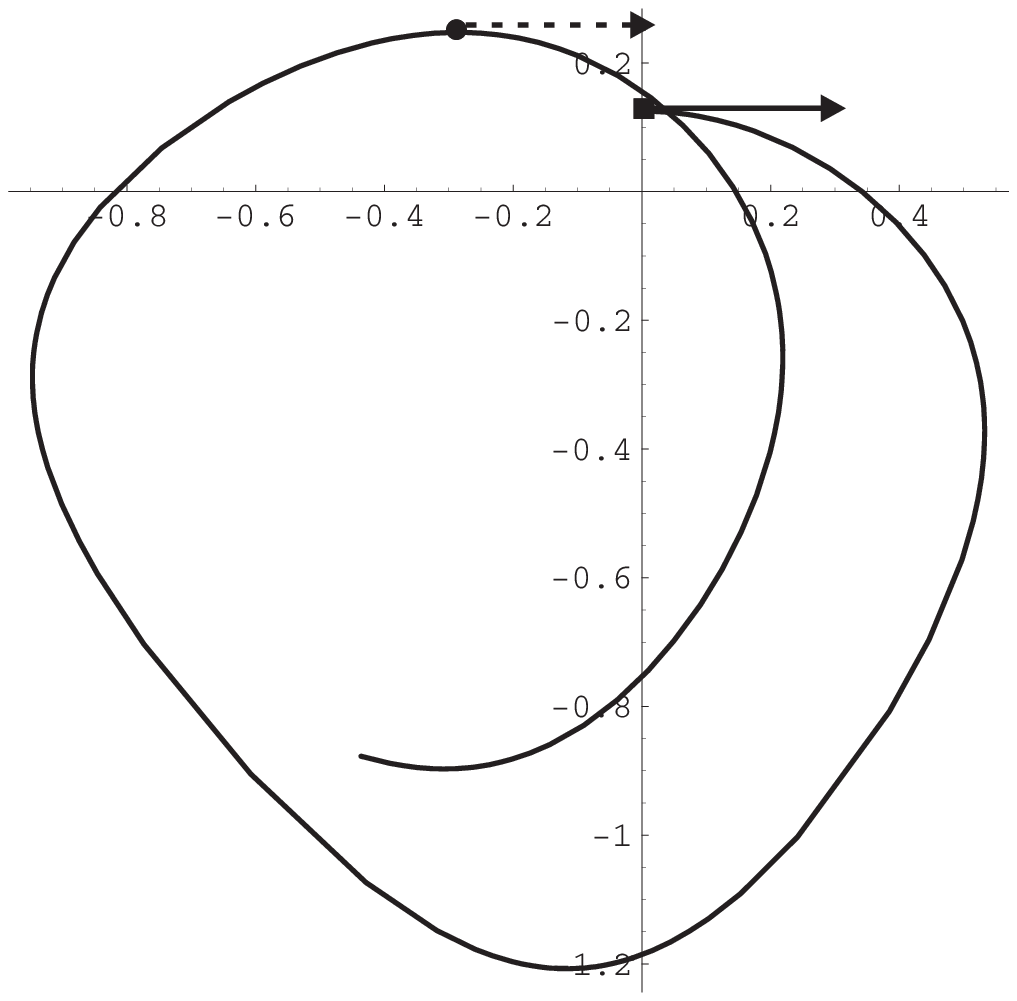}
\end{center}

\begin{center} \includegraphics[width=.15\textwidth]{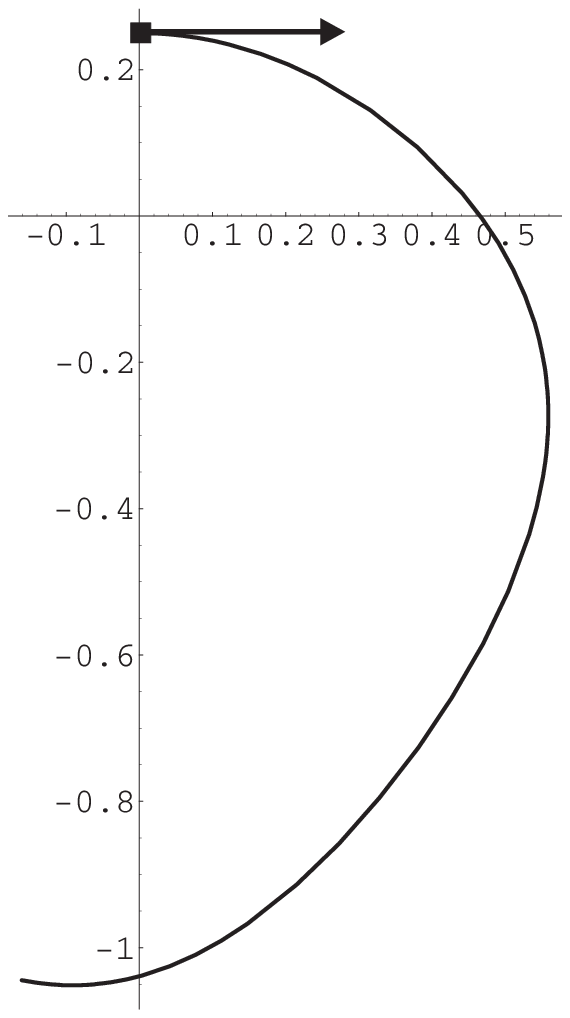}\hspace*{.5cm}\includegraphics[width=.3\textwidth]{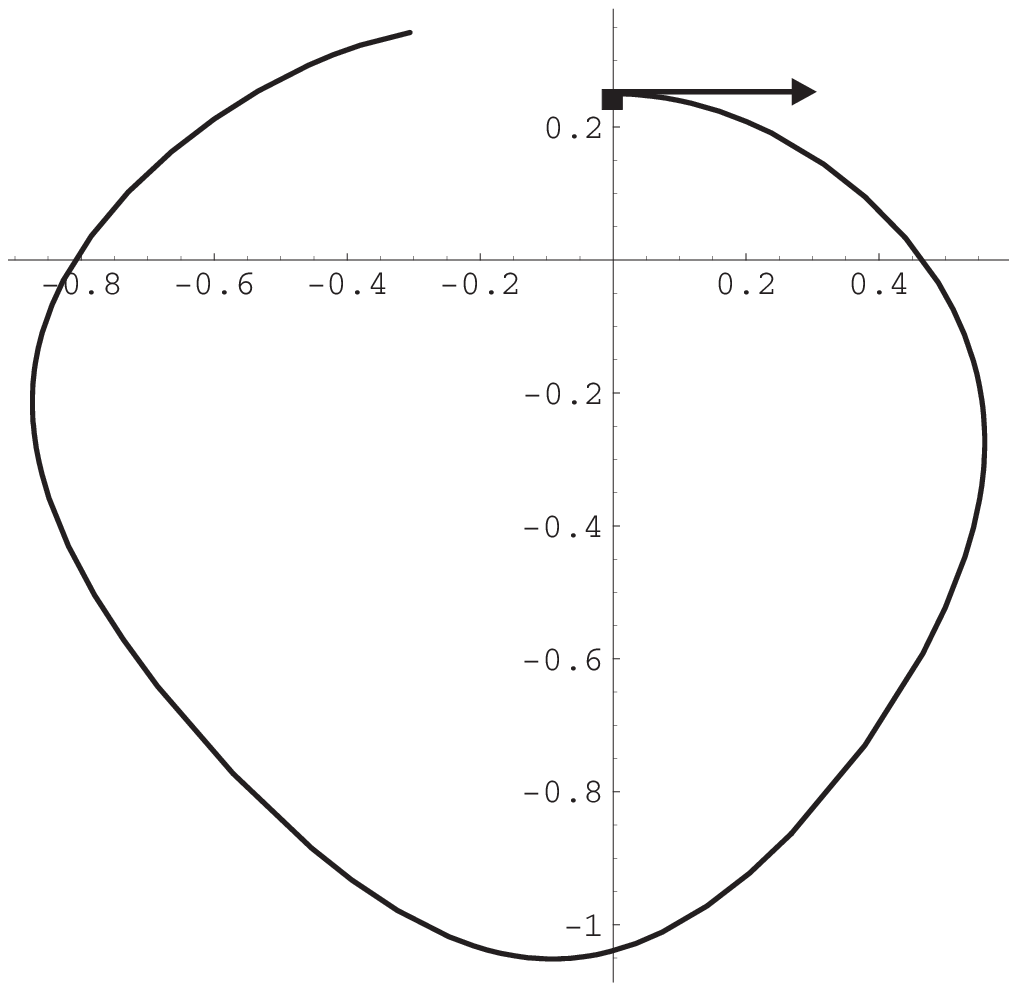}\hspace*{.5cm}
\includegraphics[width=.3\textwidth]{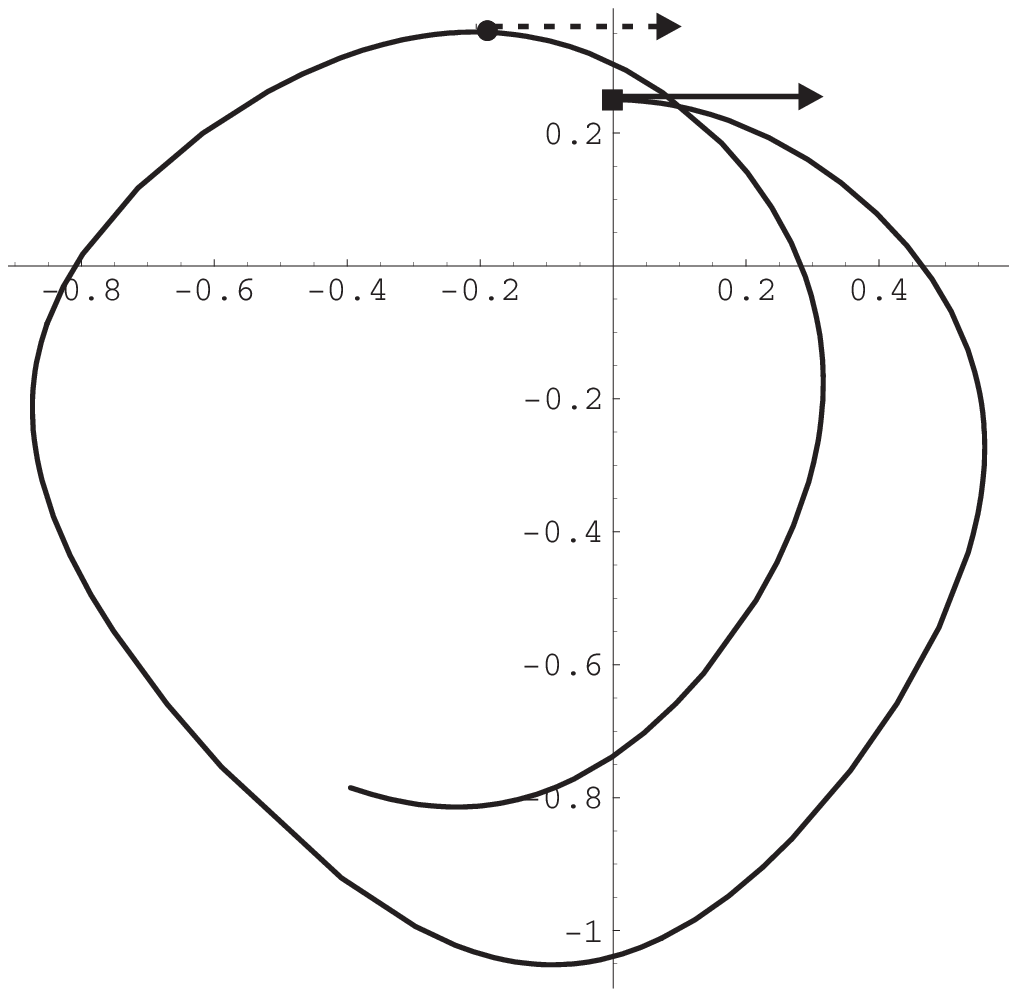}
\end{center}

\begin{center} \includegraphics[width=.15\textwidth]{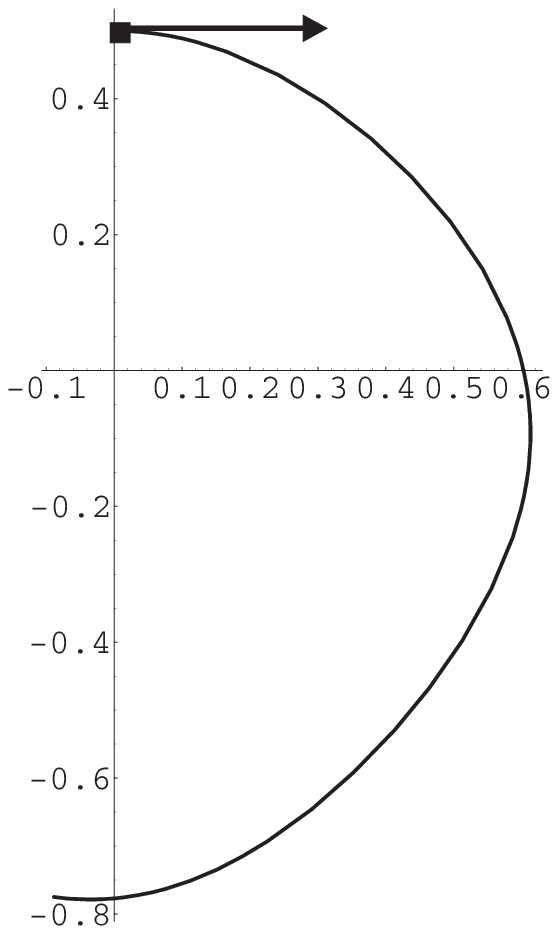}\hspace*{.5cm}\includegraphics[width=.3\textwidth]{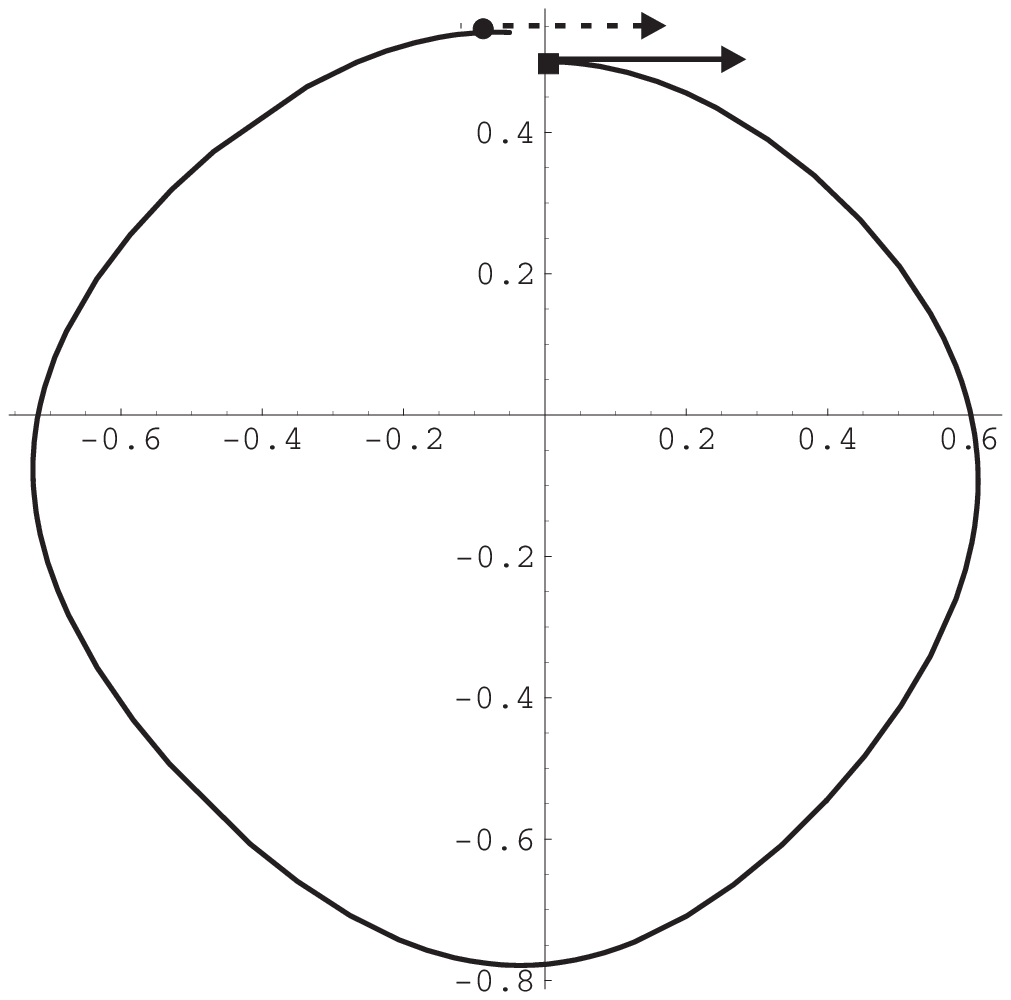}\hspace*{.5cm}
\includegraphics[width=.3\textwidth]{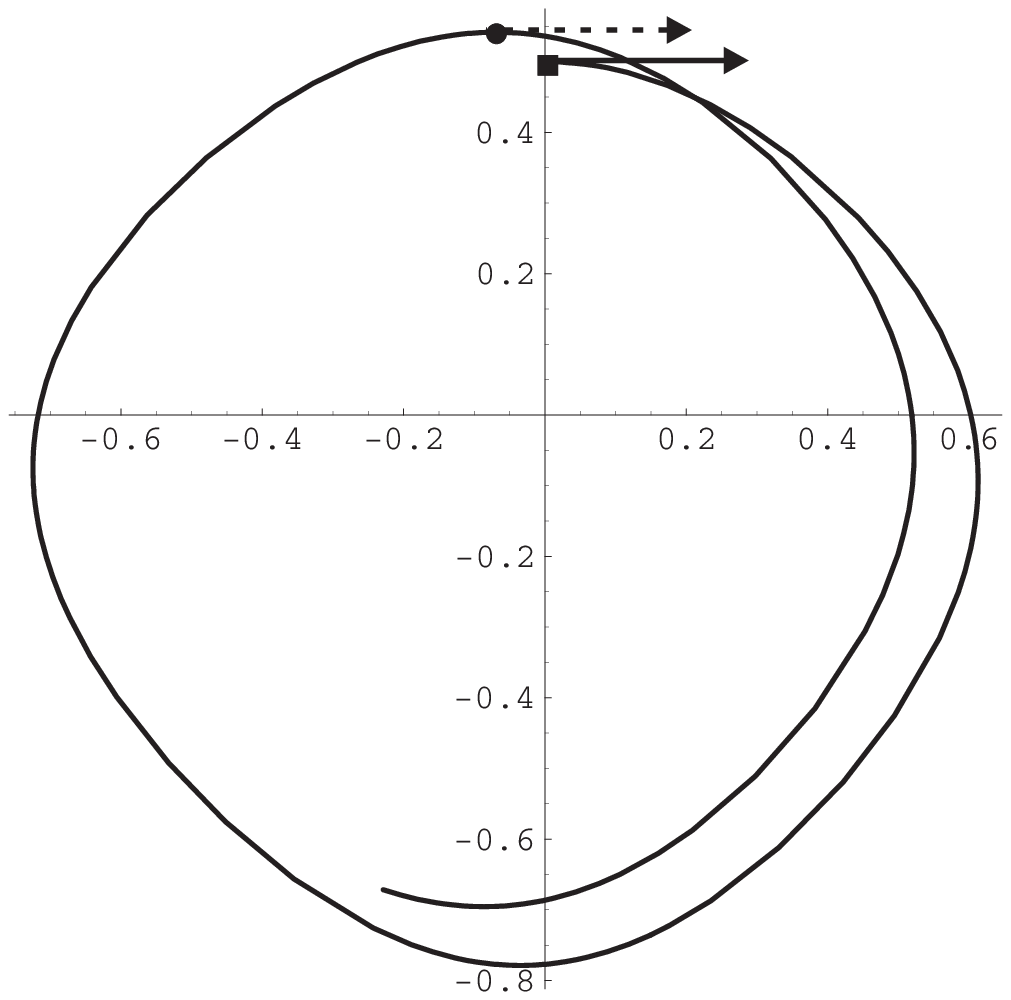}
\end{center}

\begin{center} \includegraphics[width=.15\textwidth]{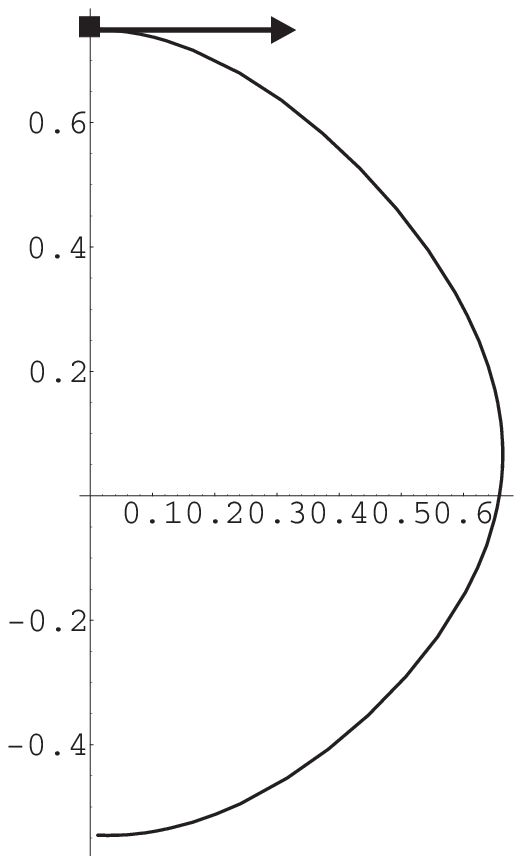}\hspace*{.5cm}\includegraphics[width=.3\textwidth]{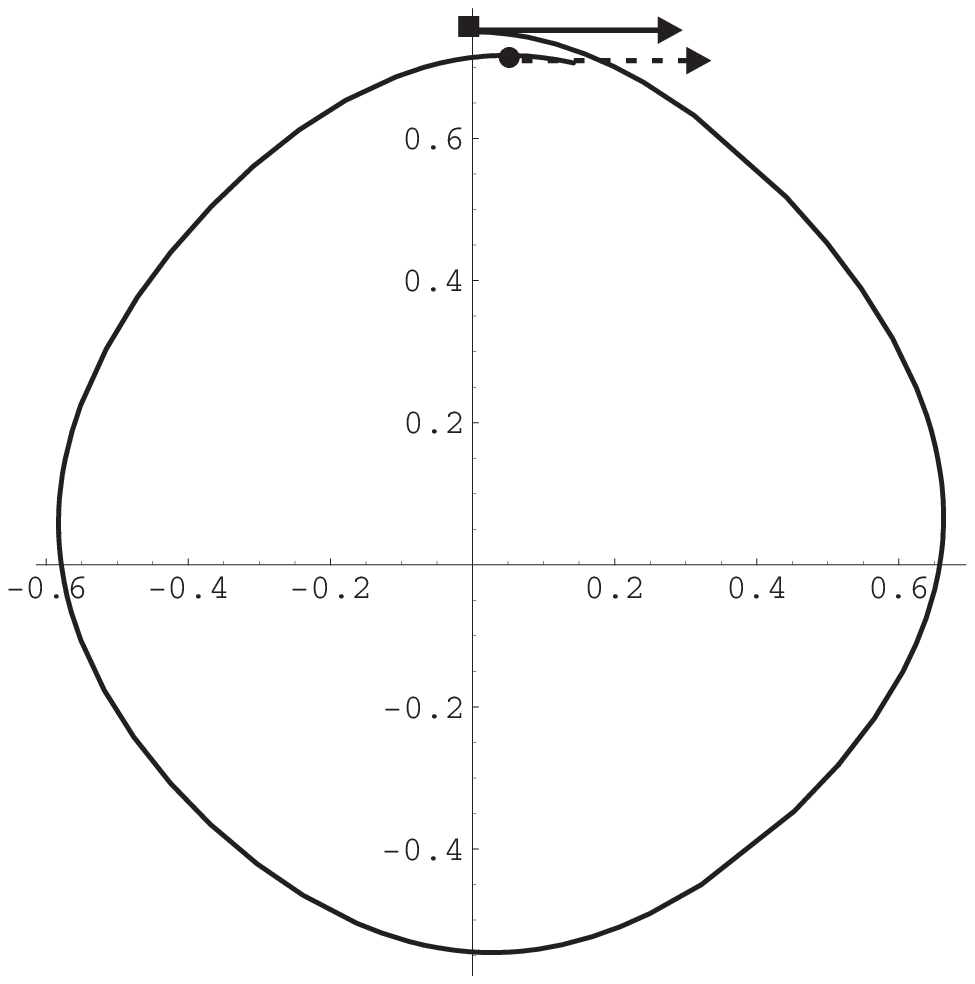}\hspace*{.5cm}
\includegraphics[width=.3\textwidth]{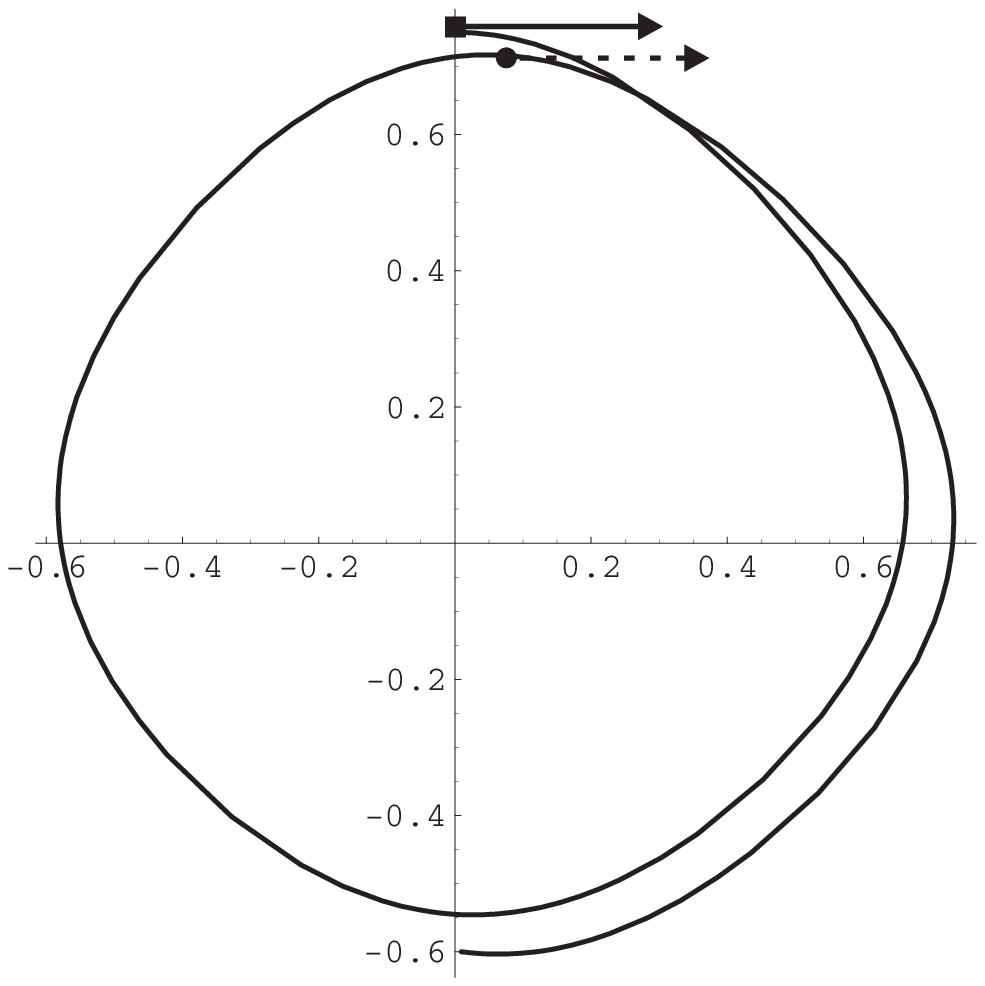}
\end{center}

{\bf Experimental result 3}. {\it Let $H\not=0$. Then there exists a simply closed curve in the $xy$-plane which is the generating curve of a $3$-invariant surface with constant mean curvature.
}

\begin{figure}[hbtp]
\begin{center} \includegraphics[width=.3\textwidth]{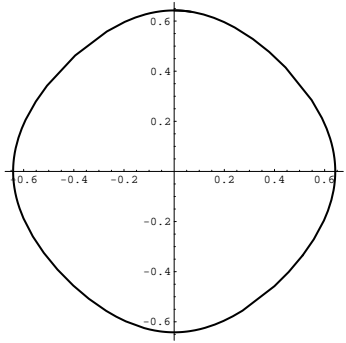}
\end{center}
\caption{The generating curve of a $3$-invariant surface with $H=1$. The curve is closed. Here the initial conditions are $x(0)=0$ and $y(0)\simeq 0.6425$. We remark  the symmetry of $\alpha$ with respect to straight-line $y=\pm x$.}\label{f-h2}
\end{figure}

 Although we do not know a proof of this statement, we have carried out different numerical graphics that seem to indicate that this result is true. Moreover, the pictures obtained of the closed solution suggest that the curve is also symmetric with respect to the straight-lines $y=\pm x$. For example, in Fig. \ref{f-h2}, the intersection points with the coordinate points are equidistant from the origin. Recall that at the origin of Sol$_3$, the directions of both straight-lines together the vertical direction $(0,0,1)$ correspond with the directions where the Ricci tensor diagonalizes.

\begin{figure}[hbtp]
\begin{center} \includegraphics[width=.4\textwidth]{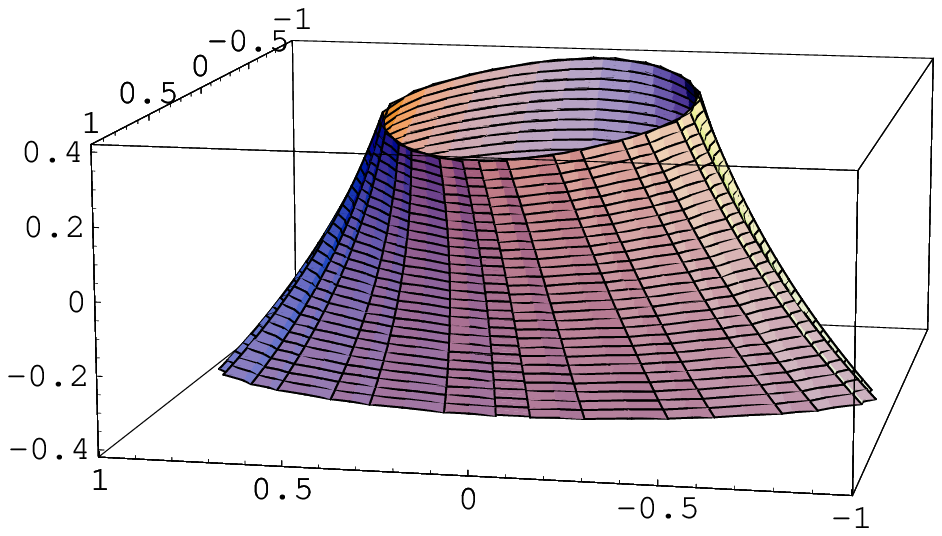}\hspace*{.5cm}\includegraphics[width=.4\textwidth]{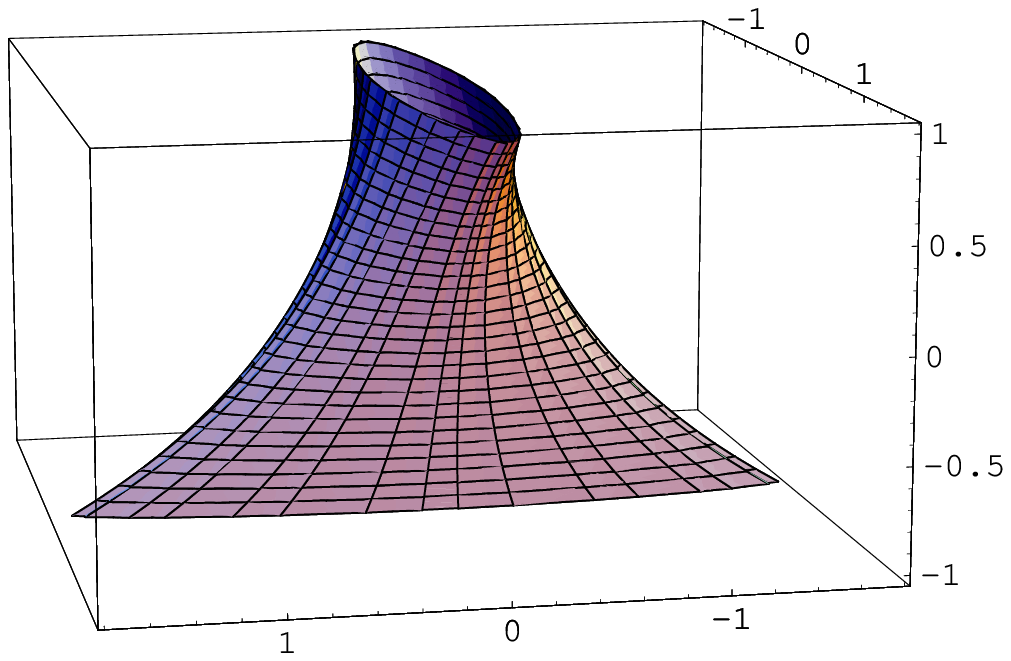}
\end{center}
\caption{Different views of the $3$-invariant surface obtained by the generating curve of Fig. \ref{f-h2}.}
\end{figure}

\section{Surfaces with constant Gauss curvature}\label{s-4}

In this section we present an example of a $3$-invariant surface in Sol$_3$ with zero Gauss curvature. Surfaces in the homogenous space  Sol$_3$ with constant Gauss curvature have not received yet the interest of geometers. In this sense, the surfaces  with constant Gauss curvature $K$   which are invariant by the groups ${\cal G}_1$ and
${\cal G}_2$ were studied in \cite{lm2} and the corresponding surfaces with $K=0$ were all obtained. The invariance of the surface makes that the Gauss curvature equation converts into an ordinary differential equation which can be partially studied.

In our context of $3$-invariant surfaces, we compute the Gauss curvature $K$. The expression of $K$ is
\begin{equation}\label{gauss}
K=K_{ext}+K(\psi_s\wedge \psi_t),
\end{equation}
where $K_{ext}$ is the extrinsic curvature of the surface and   $K(\psi_s\wedge \psi_t)$ is the sectional curvature of the tangent plane at $\psi(s,t)$. The value $K_{ext}$ is
$$K_{ext}=\frac{eg-f^2}{EG-F^2}$$
and the sectional curvature $K(\psi_s\wedge \psi_t)$ is given by
$$K(\psi_s\wedge \psi_t)=\frac{\langle\overset{\sim}{\nabla}_{\psi_s}\overset{\sim}{\nabla}_{\psi_t} \psi_t-\overset{\sim}{\nabla}_{\psi_t}\overset{\sim}{\nabla}_{\psi_s} \psi_t-\overset{\sim}{\nabla}_{[\psi_s,\psi_t]}\psi_t,\psi_s\rangle}{W}.$$
From the computations of the first and second fundamental form realized in Sect. \ref{s-2}, we have
\begin{equation}\label{ext}
K_{ext}=-\frac{A^2(x\cos\theta+y\sin\theta)^2+(\theta'+A\cos(2\theta))
\left(-x\sin\theta+y\cos\theta+A(-x^2+y^2)\right)}{W^2}
\end{equation}
On the other hand,
\begin{eqnarray*}
& &  \langle\overset{\sim}{\nabla}_{\psi_s}\overset{\sim}{\nabla}_{\psi_t} \psi_t,\psi_s\rangle =-1+\cos(2\theta)(y^2-x^2)\\
& &  \langle\overset{\sim}{\nabla}_{\psi_t}\overset{\sim}{\nabla}_{\psi_s} \psi_t,\psi_s\rangle =-(x\cos\theta+y\sin\theta)^2\\
& &  [ \psi_s,\psi_t ]  =  0.
  \end{eqnarray*}
Hence we have
\begin{equation}\label{f-we}
K(\psi_s\wedge \psi_t)=\frac{-1+A^2}{W}.
\end{equation}
Thus, by combining \eqref{ext} and \eqref{f-we},   we obtain from \eqref{gauss}
\begin{equation}\label{gauss2}
K=\frac{-(\theta'+A\cos(2\theta))
\left(-x\sin\theta+y\cos\theta+A(-x^2+y^2)\right)-1-\cos(2\theta)(x^2-y^2)A^2}{W^2}.
\end{equation}
 We search $3$-invariant surfaces in Sol$_3$ with constant Gauss curvature. Assume that $\theta$ is a constant function, $\theta(s)=\theta_0$. When the value of $\theta_0$ is $0$ or $\pi/2$, we obtain from the expression of the Gauss curvature:
\begin{enumerate}
\item If $\theta_0=0$,   we get  $K=-\frac{1}{1+y_0^2}$.  Thus the surface has constant negative Gauss curvature. Because $\theta(s)=0$, the generating curve is a straight-line parallel to the $x$-axis.  In the particular case that $y_0=0$, the surface is the hyperbolic plane $P$, which it has Gauss curvature $K=-1$, as it is well known.
\item If $\theta_0=\pi/2$, the result is similar, but now $K=-\frac{1}{1+x_0^2}$. Now the generating curves are
straight-lines parallel to the $y$-axis and when $x_0=0$, the surface is the plane $Q$.
    \end{enumerate}
These surfaces are the surfaces of type I and II that appeared in Prop. \ref{pr33}, which we know that they  are also invariant by the uniparametric groups ${\cal G}_1$ and ${\cal G}_2$. Finally, these surfaces appeared in the classification given in \cite{lm2} of constant Gauss curvature surfaces invariant by both groups.

\begin{figure}[hbtp]
\begin{center} \includegraphics[width=.8\textwidth]{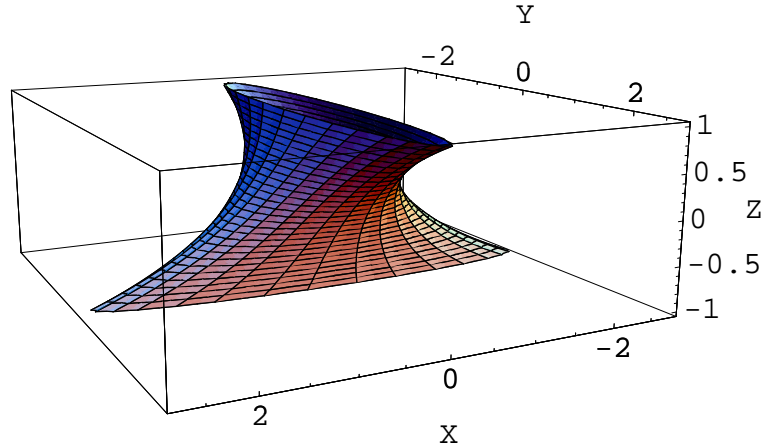}
\end{center}
\caption{A $3$-invariant surface with zero Gauss curvature, which it is generated by the curve
$\alpha(s)=(\sin(s),-\cos(s),0)$.}\label{k0}
\end{figure}

More surprisingly is that we can obtain a surface with zero Gauss curvature. From \eqref{gauss2}, if $K=0$,   the generating curve satisfies \eqref{e1} and
\begin{equation}\label{e-k0}
(\theta'+A\cos(2\theta)) \left(-x\sin\theta+y\cos\theta+A(-x^2+y^2)\right) +1+\cos(2\theta)(x^2-y^2)A^2=0.
 \end{equation}
Although Eq. \eqref{e-k0} is a bit cumbersome, a solution of this equations is given by an Euclidean circle in the $xy$-plane centered at the origin. Exactly,   if $r>0$, take
$$x(s)=r\sin(\frac{s}{r}),\ y(s)=-r\cos(\frac{s}{r}), \ \theta(s)=\frac{s}{r}.$$
Then it is direct that $\{x(s),y(s),\theta(s)\}$ is a solution of \eqref{e1}-\eqref{e-k0} with initial conditions $x(0)=0$, $y(0)=-r$ and $\theta(0)=0$. A picture of the corresponding surface appears in  Fig. \ref{k0}.

\emph{Acknowledgements}: The author thanks the valuable  discussions with Joaqu\'{\i}n P\'erez.  Part of this work was realized while the author was visiting  the Department of Mathematics of the Universidade Federal do Ceara at Fortaleza in September of 2011, whose hospitality is gratefully acknowledged.


\end{document}